\documentclass[a4paper,10pt]{amsproc}

\usepackage{amsfonts, amsmath, amssymb, graphicx, mathrsfs}
\usepackage[all]{xy}
\usepackage{tikz-cd}
\newdir{ >}{!/6pt/@{ }*:(1,0.2)@_{>}*:(1,-0.2)@^{>}}

\theoremstyle{plain}

\newtheorem{theorem}{Theorem}[section]
\newtheorem{lemma}[theorem]{Lemma}
\newtheorem{proposition}[theorem]{Proposition}
\theoremstyle{definition}
\newtheorem{definition}[theorem]{Definition}

\newtheorem{remark}[theorem]{Remark}
\newtheorem{counter example}[theorem]{Counter Example}

\newtheorem{corollary}[theorem]{Corollary}
\newtheorem{example}[theorem]{Example}

\numberwithin{equation}{section}

\begin{document}

\title[Concerning semirings of measurable functions ]{Concerning semirings of measurable functions}

\author[P. Biswas]{Pronay Biswas}
\address{Department of Mathematics, Jadavpur University, 188, Raja Subodh Chandra Mallick Rd, Kolkata 700032, West Bengal}
\email{pronayb.math.rs@jadavpuruniversity.in, pronaybiswas1729@gmail.com}
\author[S. Bag]{Sagarmoy Bag}
\address{Department of Mathematics, Bankim Sardar College, Tangrakhali, West Bengal 743329, India}
\email{sagarmoy.bag01@gmail.com}
\author[S. K. Sardar]{Sujit Kumar Sardar*\footnote{ * Corresponding author}}
\address{Department of Mathematics, Jadavpur University, 188, Raja Subodh Chandra Mallick Rd, Kolkata 700032, West Bengal}
\email{sksardarjumath@gmail.com, sujitk.sardar@jadavpuruniversity.in}

\subjclass[2010]{Primary 54C40; Secondary 46E30}

\large
\keywords{$z$-congruence, maximal congruence, prime $z$-congruence, compact measurable space, real maximal congruence}
\thanks {}

\maketitle

\begin{abstract}
		For a measurable space $(X,\mathcal{A})$, let $\mathcal{M}^+(X,\mathcal{A})$ be the commutative semiring of non-negative real-valued measurable functions with pointwise addition and pointwise multiplication. We show that there is a lattice isomorphism between the ideal lattice of $\mathcal{M}^+(X,\mathcal{A})$ and the ideal lattice of its ring of differences $\mathcal{M}(X,\mathcal{A})$.
  Moreover, we infer that each ideal of $\mathcal{M}^+(X,\mathcal{A})$ is a semiring $z$-ideal. 
  We investigate the duality between cancellative congruences on $\mathcal{M}^{+}(X,\mathcal{A})$ and $Z_{\mathcal{A}}$-filters on $X$. We observe that for $\sigma$-algebras, compactness and pseudocompactness coincide, and we provide a new characterization for compact measurable spaces via algebraic properties of $\mathcal{M}^+(X,\mathcal{A})$. It is shown that the space of (real) maximal congruences on $\mathcal{M}^+(X,\mathcal{A})$ is homeomorphic to the space of (real) maximal ideals of the $\mathcal{M}(X,\mathcal{A})$. We solve the isomorphism problem for the semirings of the form $\mathcal{M}^+(X,\mathcal{A})$ for compact and realcompact measurable spaces. 
	\end{abstract}

	
	\section{Introduction}

        In what follows, the pair $(X,\mathcal{A})$ stands for a nonempty set $X$ with a $\sigma$-algebra $\mathcal{A}$ on $X$. We call $(X,\mathcal{A})$ a measurable space. A $\sigma$-algebra $\mathcal{A}$ is said to separate points if for any two distinct points $x, y\in X$, we get $A\in \mathcal{A}$ such that $x\in A$ and $y\notin A$. Unless otherwise stated, by a measurable space we shall always mean a \emph{$T$-measurable space} (\cite{Estaji}), that is, $\mathcal{A}$ separates points of $X$. A function $f\colon X\rightarrow \mathbb{R}$ is said to be $\mathcal{A}$-measurable (or measurable) if $f^{-1}(\mathfrak{O})\in \mathcal{A}$, where $\mathfrak{O}$ is any open set in $\mathbb{R}$. The collection of all real-valued measurable functions on $(X,\mathcal{A})$, denoted by $\mathcal{M}(X,\mathcal{A})$, with pointwise addition and pointwise multiplication, forms a commutative lattice-ordered ring with unity. 

        In this paper, we initiate a study of the \emph{positive cone} (the set of all non-negative elements) of the ring $\mathcal{M}(X,\mathcal{A})$, which we denote by $\mathcal{M}^+(X,\mathcal{A})$. The set $\mathcal{M}^+(X,\mathcal{A})$ forms a commutative lattice-ordered semiring with the usual operations. One of the main themes of this paper is to constitute various bridges between the ideals and congruences of the ring $\mathcal{M}(X,\mathcal{A})$ and the semiring $\mathcal{M}^+(X,\mathcal{A})$. 

        In some recent papers like \cite{AA2020, A2020, Estaji}, the ring $\mathcal{M}(X,\mathcal{A})$ has been studied extensively. It is easy to show that $\mathcal{M}(X,\mathcal{A})$ is always a von Neumann regular ring. Therefore each ideal of $\mathcal{M}(X, \mathcal{A})$ is a $z$-ideal in the sense of Mason \cite{Mason}. Estaji et. al. \cite{Estaji} gave a complete description of the maximal ideals in terms of the lattice-theoretic aspects of $\mathcal{A}$. They have solved the isomorphism problem of rings like $\mathcal{M}(X,\mathcal{A})$ in the category of compact measurable spaces.
        Acharyya et. al. \cite{A2020} yielded an alternative method to describe the maximal ideals via the space of all $Z_\mathcal{A}$-ultrafilters on $X$. They showed that the structure space of the ring $\mathcal{M}(X,\mathcal{A})$ is zero-dimensional. 
        In \cite{AA2020}, a proof of the isomorphism problem for rings like $\mathcal{M}(X,\mathcal{A})$ is given, in the category of realcompact measurable spaces. In this paper, we try to embark on an alternative study of compact measurable spaces and realcompact measurable spaces in a more congruence-theoretic slant.
	
	In section \ref{P3, Section 3}, we provide a complete description of ideals of the semiring $\mathcal{M}^+(X,\mathcal{A})$. A major achievement of this section is an isomorphism between the ideal lattices of the semiring $\mathcal{M}^+(X,\mathcal{A})$ and the ring $\mathcal{M}(X,\mathcal{A})$. Whence each ideal is of the form $I\cap \mathcal{M}^+(X,\mathcal{A})$, for some ideal $I$ of the ring $\mathcal{M}(X,\mathcal{A})$. Under this circumstance, each ideal of $\mathcal{M}^+(X,\mathcal{A})$ turns out to be a $z$-ideal in the sense of \cite{Biswas2024}. 

        In section \ref{P3, Section 4} we study the interplay between cancellative congruences on $\mathcal{M}^+(X,\mathcal{A})$ and $Z_\mathcal{A}$-filters on $X$. We define $z$-congruences on both $\mathcal{M}^+(X,\mathcal{A})$ and $\mathcal{M}(X,\mathcal{A})$, which are heavily related to the concept of zero-sets. As anticipated, it turns out that in the ring $\mathcal{M}(X,\mathcal{A})$, there is a one-one correspondence between $z$-ideals and $z$-congruences. Consequently, we give a purely algebraic description of $z$-congruences on $\mathcal{M}^+(X,\mathcal{A})$. 

        In section \ref{P3, Section 5}, by exploiting the duality of maximal congruences and $Z_{\mathcal{A}}$-ultrafilters we show that the structure space of the semiring $\mathcal{M}^+(X,\mathcal{A})$ is homeomorphic to the structure space of the ring $\mathcal{M}(X,\mathcal{A})$. We achieve another characterization of compact measurable spaces via the semiring $\mathcal{M}^+(X,\mathcal{A})$. An interesting fact here is that in the case of $\sigma$-algebras (viz. $\sigma$-frames), the concepts of compactness and pseudocompactness coincide. 

        In section \ref{P3, Section 6} our purpose is twofold. First, we initiate a study on quotients of the semiring $\mathcal{M}^+(X,\mathcal{A})$. We show that $\mathcal{M}^+(X,\mathcal{A})/\rho$ is a totally ordered semiring if $\rho$ is a maximal congruence and the quotient semiring  $\mathcal{M}^+(X,\mathcal{A})/\rho$ is either isomorphic to or, it properly contains the semifield of non-negative reals. This leads us to the definition of real maximal congruences. We observe that the collection of all real maximal congruences, denoted by $\mathcal{RM}Cong(\mathcal{M}^+(X,\mathcal{A}))$, can be perceived dually as a topological space with the Stone topology and as $T$-measurable space. In both cases $\mathcal{RM}Cong(\mathcal{M}^+(X,\mathcal{A}))$ is homeomorphic to $RMax(X,\mathcal{A})$ (set of all real maximal ideals of $\mathcal{M}(X,\mathcal{A})$) as a topological space and as a measurable space. Lastly, we solve the isomorphism problem for the semiring of the form $\mathcal{M}^+(X,\mathcal{A})$ in the category of realcompact measurable spaces.

	\section{Preliminaries}
	
	To make this article self-contained, we recall some basics from semiring theory.
	
	A \textit{semiring} $S$ is a non-empty set with two binary operations $+$ and $\cdot$ such that $(S,+)$ and $(S,\cdot)$ are commutative monoids and $(a+b)\cdot c=a\cdot c+b\cdot c$ and $a\cdot 0=0$, for all $a,b,c\in S$.
	
	A semiring $S$ is said to be an \textit{additively cancellative semiring} (or simply \emph{cancellative}) if $a+c=b+c$ implies $a=b$, for all $a,b,c\in S$. 

        For a cancellative semiring $(S,+\cdot)$, we define $D(S)=\{a-b\colon a,b\in S\}$. Then $(D(S),+,\cdot)$ forms a ring containing the formal differences of elements from $S$. We call $D(S)$ the \emph{ring of differences} of the cancellative semiring $S$ (cf. Theorem 5.11, \cite{Hebisch}).
	
	An \textit{ideal} $I$ of $S$ is a submonoid of $(S,+)$ such that $s\cdot t\in I$, for all $s\in S$ and for all $t\in I$. 
	
	Semiring, being a more general algebraic structure than a ring, contains more classes of ideals than a ring. An ideal $I$ is said to be a \textit{$k$-ideal} if $a+b\in I$ and $b\in I$ implies $a\in I$. The class of $k$-ideals behaves more like ring ideals. Lastly, we call an ideal $I$ a \textit{strong ideal} if $a+b\in I$, then both $a\in I$ and $b\in I$.

	Unlike rings, the factor objects of semirings are not determined by ideals. Instead of ideals, congruence plays an important role in the quotient. 
	
	\begin{definition}
		A \textit{congruence} $k$ is an equivalence relation on $S$, which is also a subsemiring of the product semiring $S\times S$.  
	\end{definition}
	Equivalently, a congruence is an equivalence relation on $S$ which is compatible with the binary operations. By compatibility, we mean:
	\begin{center}
	\begin{itemize}
		\item[1.] $(a,b)\in k$ and $(c,d)\in k\implies$ $(a+c,b+d)\in k$.
		\item[2.] $(a,b)\in k$ and $(c,d)\in k \implies (a\cdot c, b\cdot d)\in k$.
	\end{itemize}
	\end{center}
	A congruence $k$ on $S$ is said to be a \textit{cancellative congruence} $(a+c,b+c)\in k$ implies $(a,b)\in k$, for all $a,b,c\in S$.
	
	A cancellative congruence $\rho$ on $S$ is said to be a \textit{regular congruence} if there exist elements $e_1$ and $e_2$ in $S$ such that for all $a\in S$, $(a+e_1a, e_2a)\in \rho$ and $(a+ae_1, ae_2)\in \rho$.
	
	Evidently, if $S$ is a commutative semiring with unity then the class of all regular congruences coincides with the class of all cancellative congruences. 
	
	\begin{definition}
		A semiring $S$ with a partial order $'\leq '$ is called a partially ordered semiring if the following conditions are satisfied: for all $a,b,c,d\in S$,
		
		\begin{itemize}
			\item[1)] $a\leq b \Leftrightarrow a+c\leq b+c$.
			
			\item[2)] $a\leq c, b\leq d\Rightarrow ad+cd\leq ab+cd$.
		\end{itemize}
	\end{definition}
	
	\begin{definition}
		A congruence $\rho$ on a partially ordered semiring $S$ is called convex if for all $a, b, c,d\in S, (a,b)\in \rho, a\leq c\leq d\leq b\Rightarrow (c,d)\in \rho$.
	\end{definition}
	
	
	
	The following theorem is noted in \cite{A1995}.
	
	\begin{theorem}\label{P3, Th2.4}
		Let $S$ be a partially ordered semiring and $\rho$ be a regular congruence on $S$. Then $S/\rho$ is a partially ordered semiring, according to the definition $\rho (a)\leq \rho (b)$ if and only if there exists $x, y\in S$ such that $(x,y)\in \rho$ and $a+x\leq b+y$, it is necessary and sufficient that $\rho$ is convex.
	\end{theorem}
	
	For any two elements $(x_1,x_2)$ and $(y_1,y_2)$ of the semiring $S\times S$, we define the \textit{twisted product} $(x_1,x_2)\cdot_t (y_1,y_2)$ as follows:
	
	\begin{center}
		$(x_1,x_2)\cdot_t (y_1,y_2)=(x_1y_1+x_2y_2, x_1y_2+x_2y_1)$
	\end{center}
	
	\begin{definition}
		A congruence $\rho$ on $S$ is called \textit{prime congruence} if for all $a, b, c, d\in H$, $(x_1,x_2)\cdot_t (y_1,y_2)\in \rho$ implies $(x_1,x_2)\in \rho$ or $(y_1,y_2)\in \rho$.
	\end{definition}
	
	The family $Cong(S)$ of all congruences on a semiring $S$ forms a complete lattice with the following operations:
	
	\begin{itemize}
		\item[1.] For any nonempty family $\mathfrak{F}$ of congruences on $S$, $\wedge \mathfrak{F}$ is defined by $(a,b)\in \wedge \mathfrak{F}$ if and only if $(a,b)\in \rho$ for every $\rho$ in $\mathfrak{F}$.
		\item[2.]  For any nonempty family $\mathfrak{F}$ of congruences on $S$, $\vee \mathfrak{F}$ is defined by $(a,b)\in \vee \mathfrak{F}$ if and only if there exist elements $a=c_0,c_1,\cdots, c_n=b$ of $S$ and congruences $\rho_1,\rho_2,\cdots, \rho_n$ of $\mathfrak{F}$ such that $(c_{i-1}, c_i)\in \rho_i$ for all $1\leq i\leq n$.

	\end{itemize}

	
	\section{The Semiring $\mathcal{M}^+(X,\mathcal{A})$ and its ideals}\label{P3, Section 3}
	
	In this section, we focus on the nature of the semiring $\mathcal{M}^+(X,\mathcal{A})$, which is the positive cone (set of all non-negative elements) of the ring $\mathcal{M}(X,\mathcal{A})$. It is easy to observe that whenever $f+g=f+h$, then $g=h$, for any $f,g,h$  in $\mathcal{M}^+(X,\mathcal{A})$. Therefore $\mathcal{M}^+(X,\mathcal{A})$ is an additively cancellative semiring. Since $\mathcal{M}^+(X,\mathcal{A})$ is the positive cone of the von Neumann regular ring $\mathcal{M}(X,\mathcal{A})$, so the semiring $\mathcal{M}^+(X,\mathcal{A})$ is a von Neumann regular semiring.

 For $f, g\in \mathcal{M}^{+}(X,\mathcal{A})$, define $f\leq g$ if and only if $f(x)\leq g(x)$ for all $x\in X$. Then $\mathcal{M}^{+}(X,\mathcal{A})$ is a partially ordered semiring with respect to the relation $'\leq'$.
	
	Also for $f,g\in \mathcal{M}^{+}(X,\mathcal{A})$, define $(f\vee g)(x)=max\{ f(x), g(x)\}$ and $(f\wedge g)(x)=min\{ f(x), g(x)\}$ for all $x\in X$. Then $f\vee g, f\wedge g\in \mathcal{M}^{+}(X,\mathcal{A})$. Therefore $\mathcal{M}^{+}(X,\mathcal{A})$ is a latticed ordered semiring with respect to $\vee$ and $\wedge$. 
	\subsection{Ideals of $\mathcal{M}^+(X,\mathcal{A})$}
	For any $f\in \mathcal{M}(X,\mathcal{A})$, define 
	\[  f^+(x)= \left\{
	\begin{array}{ll}
		f(x), & f(x)\geq 0. \\
		0, &  f(x) < 0.\\
	\end{array},
	\right.
	f^-(x)=\left\{
	\begin{array}{ll}
		0, & f(x)\geq 0. \\
		f(x), &  f(x) \leq 0.\\
	\end{array},
	\right.\]
	
	then $f=f^++f^-=f^+-(-f^-)$ and $|f|=f^+-f^-$. Clearly $f^+$ and $-f^-$ belongs to $\mathcal{M}^+(X,\mathcal{A})$. Therefore $\mathcal{M}(X,\mathcal{A})$ is the ring of differences of the semiring $\mathcal{M}^+(X,\mathcal{A})$. The following lemma manifests divisibility in the semiring $\mathcal{M}^+(X,\mathcal{A})$. 
	\begin{lemma}\label{P3, lem3.1}
		Let $f,g\in \mathcal{M}^+(X,\mathcal{A})$ and $f\leq g^r$, for some $r\geq 1$. Then $f$ is a multiple of $g$. 
	\end{lemma}
	\begin{proof}
		If $f,g\in \mathcal{M}^+(X,\mathcal{A})$ be such that $f\leq g^r$, for some $r\geq 1$. Then clearly $Z(f)\supseteq Z(g)$. Define 
		\[ h(x)=\left\{
		\begin{array}{ll}
			\frac{f(x)}{g(x)}, & x\notin Z(g). \\
			0, & x\in Z(g).\\
		\end{array},
		\right. \]
		then, both $h\vert_{Z(g)}$ and $h\vert_{X\setminus Z(g)}$ is measurable. Hence by pasting lemma, $h\in \mathcal{M}^+(X,\mathcal{A})$ and clearly $f=gh$.
	\end{proof}
	
	It is easy to observe that in a partially ordered semiring, the class of $l$-ideals coincides with the class of strong ideals. 
	
	\begin{corollary}\label{P3, cor 2.2}
		Every ideal of $\mathcal{M}^+(X,\mathcal{A})$ is a strong ideal. 
	\end{corollary}
	\begin{proof}
		Let $I$ be an ideal of $\mathcal{M}^+(X,\mathcal{A})$ and let $f+g\in I$. Then $f\leq f+g$ and $g\leq f+g$.
		By Lemma \ref{P3, lem3.1}, $f=(f+g)h$ and $g=(f+f)k$ for some $h,k\in \mathcal{M}^+(X,\mathcal{A})$. Therefore $f$ and $g$ is in $I$. Hence $I$ is a strong ideal. 
	\end{proof}
	
	\begin{remark}
		One of the contrasting natures between the semiring $\mathcal{M}^+(X,\mathcal{A})$ and the semiring $C^+(X)$ (viz. the semiring of non-negative real-valued continuous functions on a topological space $X$) is the nature of their ideals. There can exist many non-$k$-ideals in $C^+(X)$. Moreover, each ideal of $C^+(X)$ is a strong ideal (equivalently a $k$-ideal) if and only if $X$ is an $F$-space (cf. Theorem 2.1, \cite{Vechtomov}). Moreover, the lattice of ideals of $\mathcal{M}^+(X,\mathcal{A})$ is always modular (cf. Proposition 6, \cite{Alarcon}), whereas the lattice of ideals of $C^+(X)$ is modular if and only if $X$ is an $F$-space.
	\end{remark}

	By $(\mathfrak{L}(\mathcal{M}^+(X,\mathcal{A})),\vee,\wedge)$ we mean the lattice of all ideals of the semiring $\mathcal{M}^+(X,\mathcal{A})$ with $I\vee J=I+J$ and $I\wedge J=I\cap J$. Similarly $(\mathfrak{L}(\mathcal{M}(X,\mathcal{A})),\vee,\wedge)$ is the lattice of all ideals of $\mathcal{M}(X,\mathcal{A})$ with obvious join and meet. We define two maps $\alpha\colon \mathfrak{L}(\mathcal{M}(X,\mathcal{A}))\rightarrow \mathfrak{L}(\mathcal{M}^+(X,\mathcal{A}))$ and $\beta\colon \mathfrak{L}(\mathcal{M}^+(X,\mathcal{A}))\rightarrow \mathfrak{L}(\mathcal{M}(X,\mathcal{A}))$ as follows:
	
	\begin{center}
		$\alpha(I)=I\cap \mathcal{M}^+(X,\mathcal{A})$\hspace{0.5cm} and \hspace{0.5cm} $\beta(I)=\{f-g\colon f,g\in I\}$.
	\end{center}

	\begin{lemma}\label{P3, lem2.4}
		The following statements hold.
		\begin{itemize}
			\item[1.] The map $\beta$ is an onto lattice homomorphism. 
			\item[2.] The map $\alpha$ is an onto lattice homomorphism.
		\end{itemize}
	\end{lemma}
	\begin{proof}
		$1.$ The equality $\beta(I+J)=\beta(I)+\beta(J)$ and the inclusion $\beta(I\cap J)\subseteq \beta(I)\cap \beta(J)$ easily follows. Let $f\in \beta(I)\cap \beta(J)$. Then $f=g_1-h_1\in \beta(I)$ and $f=g_2-h_2\in \beta(J)$, for some $g_1,h_1\in I$. Therefore $g_1+g_2=h_1+h_2\in I\cap J$ and from Corollary \ref{P3, cor 2.2} $g_1,g_2,h_1,h_2\in I\cap J$. Therefore $f\in \beta(I\cap J)$ and $\beta(I\cap J)=\beta(I)\cap \beta(J)$. Hence $\beta$ is a lattice homomorphism. 
		
		Moreover, we show that every ideal of $\mathcal{M}(X,\mathcal{A})$ is a difference ideal. Let $f\in I$, where $I$ is an ideal in $\mathcal{M}(X,\mathcal{A})$. Clearly $f=f^+-(-f^-)$ with $f^+\leq |f|$ and $-f^-\leq |f|$. Therefore $f^+,-f^-\in I\cap \mathcal{M}^+(X,\mathcal{A})=\alpha(I)$ and $\beta(\alpha(I))=I$. Hence $\beta$ is an onto map. 
		
		$2.$ The equality $\alpha(I\cap J)=\alpha(I)\cap \alpha(J)$ and $\alpha(I)+\alpha(J)\subseteq \alpha(I+J)$ are obvious. Now suppose $f\in \alpha(I+J)$. Then $f=g+h$ and $f\leq |g|+|h|$. Applying Riesz decomposition theorem (cf. Proposition 1.1.4 of \cite{Anderson}) we get $f=s+t$, for some $s,t\in \mathcal{M}^+(X,\mathcal{A})$ such that $0\leq s\leq |g|$ and $0\leq t\leq |h|$. Recall that each ideal of $\mathcal{M}(X,\mathcal{A})$ is an $l$-ideal. Therefore $s\in I\cap \mathcal{M}(X,\mathcal{A})=\alpha(I)$ and $t\in J\cap \mathcal{M}(X,\mathcal{A})$. Hence $f\in \alpha(I)+\alpha(J)$, which validates the equality $\alpha(I+J)=\alpha(I)+\alpha(J)$. Thus we have proved that $\alpha$ is a lattice homomorphism. 
		
		In addition, for any ideal $I$ of $\mathcal{M}^+(X,\mathcal{A})$, $I\subseteq \alpha(\beta(I))$. Let $f\in \alpha(\beta(I))$. Then $f=g-h$, for some $f,g\in I$. Then $f\in I$. Indeed $I$ is a strong ideal and hence a $k$-ideal. Therefore $\alpha(\beta(I))=I$. Hence $\alpha$ is an onto map.

	\end{proof}
	
	\begin{proposition}
		The lattice $\mathfrak{L}(\mathcal{M}^+(X,\mathcal{A}))$ is isomorphic to $\mathfrak{L}(\mathcal{M}(X,\mathcal{A}))$.
	\end{proposition}
	\begin{proof}
		For any two ideals $I$ and $J$ of $\mathcal{M}^+(X,\mathcal{A})$, if $\beta(I)=\beta(J)$, then by Lemma \ref{P3, lem2.4} $I=\alpha(\beta(I))=\alpha(\beta(J))=J$. Which shows $\beta$ is injective. Similarly, $\alpha$ is also injective. 
	\end{proof}
	
	\begin{corollary}\label{P3, Cor3.6}
		Each ideal of $\mathcal{M}^+(X,\mathcal{A})$ is of the form $I\cap \mathcal{M}^+(X,\mathcal{A})$ for some ideal $I$ of $\mathcal{M}(X,\mathcal{A})$.
	\end{corollary}
	\begin{corollary}\label{P3, lem3.7}
		Each prime ideal (maximal ideal) of $\mathcal{M}^+(X,\mathcal{A})$ is of the form $P\cap \mathcal{M}^+(X,\mathcal{A})$ for some prime ideal (maximal ideal) $P$ of $\mathcal{M}(X,\mathcal{A})$.
	\end{corollary}
	
	\begin{remark}\label{P3, rem3.8}
		As a direct consequence of the above discussion and Theorem 2.11 of \cite{A2020}, we achieve a complete description of maximal ideals in $\mathcal{M}^+(X,\mathcal{A})$. Each maximal ideal of $\mathcal{M}^+(X,\mathcal{A})$ is of the form $M_+^p=\{f\in \mathcal{M}^+(X,\mathcal{A})\colon p\in cl_{\hat{X}}Z(f)\}$, where $\hat{X}$ is the space of all ultrafilters of the measure space $(X,\mathcal{A})$ under Stone-topology.
	\end{remark}
	
	\begin{definition}(\cite{Biswas2024})
		An ideal $I$ of a semiring $(S,+,\cdot,0,1)$ is said to be a $z$-ideal if $\mathcal{M}^+_a\subseteq I$ for every $a\in I$.
	\end{definition}
	Here $\mathcal{M}^+_a=\underset{M\in \mathcal{M}^+(a)}{\bigcap}M$ and $\mathcal{M}^+(a)$ is the set of all maximal ideals of $S$ containing $a$.
	
	\begin{lemma}\label{P3, lem3.10}
		For any $f\in \mathcal{M}^+(X,\mathcal{A})$, $\mathcal{M}^+_f=\{g\in \mathcal{M}^+(X,\mathcal{A})\colon Z(f)\subseteq Z(g)\}$. 
	\end{lemma}
	
	The proof relies on Lemma \ref{P3, lem3.1} and Remark \ref{P3, rem3.8}. As an easy consequence of Lemma \ref{P3, lem3.10}, we have the following.
	
	\begin{corollary}\label{P3, Cor3.11}
		Each ideal of $\mathcal{M}^+(X,\mathcal{A})$ is a $z$-ideal.
	\end{corollary}
	
	\begin{remark}
		It is evident that for any $f\in \mathcal{M}(X,\mathcal{A})$, $Z(f)=Z(|f|)$. Therefore $Z[\mathcal{M}(X,\mathcal{A})]=Z[\mathcal{M}^+(X,\mathcal{A})]$. In other words, the collections of zero-sets of the ring $\mathcal{M}(X,\mathcal{A})$ and the semiring $\mathcal{M}^+(X,\mathcal{A})$ are the same. Moreover, let $I$ be an ideal of the ring $\mathcal{M}(X,\mathcal{A})$, then $f\in I$ if and only if $|f|\in I$, so $Z[I]=Z[I\cap \mathcal{M}^+(X,\mathcal{A})]$. Now since each ideal of the semiring $\mathcal{M}^+(X,\mathcal{A})$ is of the form $I\cap \mathcal{M}^+(X,A)$, where $I$ is an ideal of the ring $\mathcal{M}(X,\mathcal{A})$ (cf. Corollary \ref{P3, Cor3.6}). We conclude that $Z[I\cap \mathcal{M}^+(X,\mathcal{A})]$ is a $Z_\mathcal{A}$-filter on $X$. In the next section, we deal with the question of whether we can extend this ideal-filter connection to a congruence-filter connection between $Cong(\mathcal{M}^+(X,\mathcal{A}))$ and $Z_\mathcal{A}$-filters on $X$. 
	\end{remark}

	
	\section{Congruences on $\mathcal{M}^+(X,\mathcal{A})$}\label{P3, Section 4}

	For any congruence $\rho$ on $\mathcal{M}^+(X,\mathcal{A})$, we define 
	\begin{center}
		$E(\rho)=\{E(f,g)\colon (f,g)\in \rho\}$,
	\end{center}
	where $E(f,g)=\{x\in X\colon f(x)=g(x)\}$, is the agreement set of $f$ and $g$. It is clear that $E(f,g)=Z(f-g)$.

	\begin{theorem}\label{P3, Th4.1}
		For a measurable space $(X,\mathcal{A})$, $A\in \mathcal{A}$ if and only if it is the agreement set of some functions $f,g$ in $\mathcal{M}^{+}(X,\mathcal{A})$.
	\end{theorem}
	
	\begin{proof}
		Let $A\in \mathcal{A}$. Then $A=Z(f)$, where $f=\chi_{A^c}\in \mathcal{M}^{+}(X,\mathcal{A})$. Therefore $A=E(f,\boldsymbol{0})$.
		
		Conversely, let $A$ be an agreement set of $f, g$ in $\mathcal{M}^{+}(X,\mathcal{A})$. Then $A=Z(f-g)\in \mathcal{A}$.
	\end{proof}
	
	Likewise, in the case of ideals, it is customary to ask questions about the structure of congruences on $\mathcal{M}^+(X,\mathcal{A})$. The following generic example shows that not all congruences on $\mathcal{M}^+(X,\mathcal{A})$ are cancellative. 
	
	\begin{example}\label{P3, exam3.14}
		Let $S=\mathcal{M}^+(X,\mathcal{A})\setminus \{\boldsymbol{0}\}$. Then $k=(S\times S)\cup \Delta_{\mathcal{M}^+(X,\mathcal{A})}$ is a non-trivial congruence on $\mathcal{M}^+(X,\mathcal{A})$, where $\Delta_{\mathcal{M}^+(X,\mathcal{A})}=\{(f,f)\colon f\in \mathcal{M}^+(X,\mathcal{A})\}$. Then $(f+\boldsymbol{1},f)\in k$ for any $f\in S$, but $(\boldsymbol{1},\boldsymbol{0})\notin k$ because zero-class of $k$ is a singleton set, that is, $[0]_{k}=\{(\boldsymbol{0},\boldsymbol{0})\}$.
	\end{example}

	An easy conclusion we can make from the above example is that, unlike the connection of ideals and $Z_\mathcal{A}$-filters, we cannot create congruences and $Z_\mathcal{A}$-filters connection. Indeed, $E(f+\boldsymbol{1},f)=\phi\in E(k)$, where $k_+$ is the congruence defined in Example \ref{P3, exam3.14}. We observe the following important correlations. Compare with the semiring $C^+(X)$; see Theorem 3.2 and Theorem 3.3 of \cite{A1993}.
	
	\begin{proposition}\label{P3, prop 4.3}The following statements hold for any measurable space $(X,\mathcal{A})$.
		\begin{itemize}
			\item[1.] If $\rho_+$ is a cancellative congruence on $\mathcal{M}^{+}(X,\mathcal{A})$, then $E(\rho_+ )=\{ E(f,g): (f,g)\in \rho_+\}$ is an $Z_\mathcal{A}$-filter on $X$.  
			\item[2.]  If $\mathfrak{F}$ is an $Z_\mathcal{A}$-filter on $X$, then $E^{-1}(\mathfrak{F})=\{ (f,g)\in \mathcal{M}^{+}(X,\mathcal{A})\times \mathcal{M}^{+}(X,\mathcal{A}): E(f,g)\in \mathfrak{F}\}$ is a cancellative congruence on $\mathcal{M}^{+}(X,\mathcal{A})$.
		\end{itemize}
	\end{proposition}

	
	\subsection{On $z$-congruences of $\mathcal{M}^+(X,\mathcal{A})$}
	
	Likewise in the case of the semiring $C^+(X)$, it is natural to consider $z$-congruences on $\mathcal{M}^+(X,\mathcal{A})$. 
	
	\begin{definition}
		A congruence $\rho$ on $\mathcal{M}^{+}(X,\mathcal{A})$ is called \emph{$z$-congruence} if for all $f,g$ in $\mathcal{M}^{+}(X,\mathcal{A})$, $E(f,g)\in E(\rho)$ implies that $(f,g)\in \rho$. 
	\end{definition}
	
	Therefore for each $z$-congruence $\rho$ on $\mathcal{M}^+(X,\mathcal{A})$, $E^{-1}(E(\rho))=\rho$. Each $z$-congruence is a cancellative congruence. The set of all $z$-congruences on $\mathcal{M}^+(X,\mathcal{A})$ is denoted by $\mathcal{Z}Cong$. Let us denote the collection of all $\mathcal{A}$-filters on $X$ by $\mathcal{Z}_\mathcal{A}$. Both $\mathcal{Z}Cong$ and $\mathcal{Z}_\mathcal{A}$ are partially ordered by inclusions. 
	
	\begin{theorem}\label{P3, Th4.5}
		The map $E\colon (\mathcal{Z}Cong,\subseteq) \rightarrow (\mathcal{Z}_\mathcal{A},\subseteq)$ is an order-isomorphism. 
	\end{theorem}
	
	The proof relies on the fact that both $E$ and $E^{-1}$ are order preserving maps and $E^{-1}(E(\rho)=\rho$ and $E(E^{-1}(\mathfrak{F}))=F$ for any $z$-congruence $\rho$ and $Z_\mathcal{A}$-filter $F$ respectively. The following class of $z$-congruences can be easily obtained by Theorem \ref{P3, Th4.5}.
	
	\begin{corollary}
		Every maximal congruence on $\mathcal{M}^+(X,\mathcal{A})$ is a $z$-congruence. 
	\end{corollary}
	
	\begin{theorem}\label{P3, Th4.6}
		An intersection of an arbitrary non-empty family of $z$-congruences on $\mathcal{M}^+(X,\mathcal{A})$ is a $z$-congruence. 
	\end{theorem}
	
	\begin{proof}
		Let $\{\rho_i\}_{i\in I}$ be a non-empty family of $z$-congruences on $\mathcal{M}^+(X,\mathcal{A})$, where $I$ is an index set. Then consider the congruence $\rho=\underset{i\in I}{\bigcap}\rho_i$ on $\mathcal{M}^+(X,\mathcal{A})$. Clearly, $E(\rho)\subset \underset{i\in I}{\bigcap} E(\rho_i)$, and hence $E^{-1}(E(\rho))\subset \underset{i\in I}{\bigcap} E^{-1}(E(\rho_i))$. Since each $\rho_i$ is a $z$-congruence, $E^{-1}(E(\rho_i))=\rho_i$ for all $i\in I$, and hence $E^{-1}(E(\rho))\subset \underset{i\in I}{\bigcap}\rho_i=\rho$. The opposite inclusion is trivial. Therefore $E^{-1}(E(\rho))=\rho$. This completes the proof. 
	\end{proof}
	
	Therefore $(\mathcal{Z}Cong,\wedge)$ can be regarded as a complete $\wedge$-semilattice, where $\rho_1\wedge \rho_2=\rho_1\cap \rho_2$. We define join of two $z$-congruences as 
	\begin{center}
		$\rho_1\vee_z \rho_2=(\rho_1\vee \rho_2)_z$, 
	\end{center}
	where $\rho_1\vee \rho_2$ is the usual join of two congruences and $(\rho_1\vee\rho_2)_z$ is the smallest $z$-congruence containing $\rho_1\vee\rho_2$. In view of Theorem \ref{P3, Th4.6} and Theorem 2.31 of \cite{Davey} we have the following result. 
	
	\begin{corollary}
		The lattice $(\mathcal{Z}Cong,\vee_z,\wedge)$ is a complete lattice. 
	\end{corollary}
	
	We have given a complete description of ideals of $\mathcal{M}^+(X,\mathcal{A})$ (cf. Corollary \ref{P3, Cor3.6}).
	Now we investigate if there is a complete description of cancellative congruences on $\mathcal{M}^+(X,\mathcal{A})$ in terms of the congruences on the ring $\mathcal{M}^+(X,\mathcal{A})$. For that we define $\mathcal{E}(f,g)=\{x\in X\colon f(x)=g(x)\}$, the agreement set of $f,g\in \mathcal{M}(X,\mathcal{A})$. There is a one-one correspondence between the set of the ideals of $\mathcal{M}(X,\mathcal{A})$ and the set of all congruences on $\mathcal{M}(X,\mathcal{A})$. Indeed, for a ring $R$, we know that there is a one-one correspondence between the class of congruences on $R$ and the class of ideals in $R$ (cf. Remark 7.6(iii), \cite{Hebisch}). Therefore for any congruence $k$ on $\mathcal{M}(X,\mathcal{A})$, $\mathcal{E}(k)=\{\mathcal{E}(f,g)\colon (f,g)\in k\}=Z[I_k]$, where $I_k=\{s-t\colon (s,t)\in k\}$ is the corresponding ideal of the congruence $k$. For any ideal $I$, the corresponding congruence is $k_I=\{(h,k)\colon h-k\in I\}$. Observe that $k_{I_k}=k$ for any congruence $k$.

	\begin{definition}
		A congruence $k$ on $\mathcal{M}(X,\mathcal{A})$ is said to be a $z$-congruence if $\mathcal{E}(f,g)\in \mathcal{E}(k)$ implies $(f,g)\in k$.
	\end{definition}
	
	\begin{theorem}\label{P3, Th4.9}
		Each congruence $k$ on $\mathcal{M}(X,\mathcal{A})$ is a $z$-congruence. 
	\end{theorem}
	\begin{proof}
		For any congruence $k$ on $\mathcal{M}(X,\mathcal{A})$ and its corresponding ideal $I_k$, we have $\mathcal{E}(k)=Z[I_k]$. Now $\mathcal{E}(f,g)=Z(f-g)\in Z[I_k]$ implies $f-g\in I_k$ because $I_k$ is a $z$-ideal. Hence $(f,g)\in k$. 
	\end{proof}
	\begin{corollary}\label{P3, Cor4.11}
		There is a one-one correspondence between $z$-ideals and $z$-congruences in the ring $\mathcal{M}(X,\mathcal{A})$.
	\end{corollary}
	Since $\mathcal{M}(X,\mathcal{A})$ is the difference ring of the cancellative semiring $\mathcal{M}^+(X,\mathcal{A})$, we can define a map $\nabla$ as follows
	
	\begin{center}
		$\nabla\colon Cong(\mathcal{M}(X,\mathcal{A}))\rightarrow Cong(\mathcal{M}^+(X,\mathcal{A}))$ 
		
		$\rho\longmapsto \rho^\nabla$
	\end{center}
	
	where $\rho^\nabla=\rho\cap (\mathcal{M}^+(X,\mathcal{A})\times \mathcal{M}^+(X,\mathcal{A}))$. Also, there exists a map $\Delta$ in the opposite direction, defined as follows
	
	\begin{center}
		$\Delta\colon Cong(\mathcal{M}^+(X,\mathcal{A}))\rightarrow Cong(\mathcal{M}(X,\mathcal{A}))$
		
		\hspace{1cm} $\rho_+\longmapsto \rho_+^\Delta$
	\end{center}
	
	where $\rho_+^\Delta=\{(f,g)\colon f-g=h-k\ \text{for some}\ (h,k)\in \rho_+\}$.
	
	\begin{proposition}\label{P3, Prop4.12}
		Each $z$-congruence of $\mathcal{M}^+(X,\mathcal{A})$ is of the form $k\cap (\mathcal{M}^+(X,\mathcal{A})\times \mathcal{M}^+(X,\mathcal{A}))$ for some congruence $k$ on $\mathcal{M}(X,\mathcal{A})$.
	\end{proposition}
	\begin{proof}
		Since every congruence $\mathcal{M}(X,\mathcal{A})$ is a $z$-congruence, it is easy to show that $\rho^\nabla$ is a $z$-congruence on $\mathcal{M}^+(X,\mathcal{A})$ for any congruence $\rho$ on $\mathcal{M}(X,\mathcal{A})$. Now suppose $\rho_+$ is a $z$-congruence on $\mathcal{M}^+(X,\mathcal{A})$. Let $\mathcal{E}(f,g)\in \mathcal{E}(\rho_+^\Delta)$. Then define $h=f-(f\wedge g)$ and $k=g-(g\wedge f)$. Clearly, $h,k\in \mathcal{M}^+(X,\mathcal{A})$. Also, we have $\mathcal{E}(f,g)=E(h,k)\in \rho_+$ and $(h,k)\in \rho_+$. But from the construction of $h$ and $k$, it is clear that $f-g=h-k$. Therefore $(f,g)\in \rho_+^\Delta$ and hence $\rho_+^\Delta$ is a $z$-congruence on $\mathcal{M}(X,\mathcal{A})$. Also, since every $z$-congruence is a cancellative congruence on $\mathcal{M}^+(X,\mathcal{A})$, we have $\rho_+^{\Delta\nabla}=\rho_+$ by Theorem 7.1 of \cite{Hebisch}. This completes the proof. 
	\end{proof}

	Due to Corollary \ref{P3, Cor3.11}, Theorem \ref{P3, Th4.9} and Corollary \ref{P3, Cor4.11}, we arrive at the following remarkable correspondence theorem for $z$-ideals and $z$-congruences in $\mathcal{M}^+(X,\mathcal{A})$. 
	\begin{proposition}\label{P3, Prop4.13}
		There is a one-one correspondence between $z$-ideals and $z$-congruences in the semiring $\mathcal{M}^+(X,\mathcal{A})$.
	\end{proposition}
	
	\begin{proof}
		Let $k$ be a $z$-congruence on $\mathcal{M}(X,\mathcal{A})$. Then the corresponding ideal (zeroth class), denoted by $I_k=\{f\colon (f,\boldsymbol{0})\in k\}$ is again a $z$-ideal of $\mathcal{M}^+(X,\mathcal{A})$. We denote this map from $\mathcal{Z}Cong(\mathcal{M}^+(X,\mathcal{A}))$ to $\mathfrak{L}(\mathcal{M}^+(X,\mathcal{A}))$ by $\mathcal{I}_+$. For any $z$-congruence $k$ we can easily observe that $\mathcal{I}_+(k)=\alpha(\mathcal{I}(k^\Delta))$.
		The following diagram captures the essence of our goal.
    \begin{center}



			\includegraphics{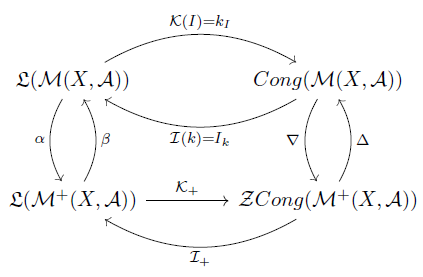}
\end{center}

Now we define a map $\mathcal{K}_+\colon \mathfrak{L}(\mathcal{M}^+(X,\mathcal{A}))\mapsto \mathcal{Z}Cong(\mathcal{M}^+(X,\mathcal{A}))$ by $\mathcal{K}_+(J)=\mathcal{K}_{\beta(I)}^\nabla$. Since each ideal $J$ of $\mathcal{M}^+(X,\mathcal{A})$ is a strong ideal, we immediately have $\mathcal{K}_{\beta(I)}^\nabla=\{(f,g)\colon f+s=g+t,\ \text{for some}\ s,t\in I\}$, which is again a $z$-congruence on $\mathcal{M}^+(X,\mathcal{A})$ corresponding to the $z$-ideal $I$. Moreover $\mathcal{I}_+(\mathcal{K}_+(I))=I$ for any ideal $I$ and $\mathcal{K}_+(I_+(k))=k$ for any $z$-congruence $k$ on $\mathcal{M}^+(X,\mathcal{A})$. This completes the proof.

	\end{proof}

	\subsection{Prime $z$-congruences}\hfill

	\begin{definition}
		An $Z_\mathcal{A}$-filter $\mathfrak{F}$ on $X$ is said to be prime if $A, B\in \mathcal{A}$, $A\cup B\in \mathfrak{F}$ implies $A\in \mathfrak{F}$ or $B\in \mathfrak{F}$.
	\end{definition}
	
	\begin{lemma}\label{P3. lemma 4.13}
		
		For any $f_1,f_2,g_1$ and $g_2$ in $\mathcal{M}^+(X,\mathcal{A})$, $E(f_1,g_1)\cup E(f_2,g_2)=E(f_1f_2+g_1g_2, f_1g_2+f_2g_1)$.
		
	\end{lemma}
	\begin{theorem}\label{P3, Th 4.14}
		If $\rho$ is a prime $z$-congruence on $\mathcal{M}^{+}(X,\mathcal{A})$, then $E(\rho)$ is a prime $z$-filter on $X$.
	\end{theorem}
	
	\begin{proof}
		Let $\rho$ be a prime $z$-congruence on $\mathcal{M}^{+}(X,\mathcal{A})$. Then by Proposition \ref{P3, prop 4.3}, $E(\rho)$ is a $Z_\mathcal{A}$-filter on $X$. Let $A_1\cup A_2\in E(\rho)$, where $A_1, A_2\in \mathcal{A}$. Then by Theorem \ref{P3, Th4.1} $A=E(f_1,g_1), A_2=E(f_2,g_2)$, for some $f_1, f_2, g_1, g_2\in \mathcal{M}^{+}(X,\mathcal{A})$. Then $A_1\cup A_2=E(f_1,g_1)\cup E(f_2,g_2)=E(f_1f_2+g_1g_2, f_1g_2+f_2g_1)\in E(\rho)$. Thus $(f_1f_2+g_1g_2, f_1g_2+f_2g_1)\in \rho$ as $\rho$ is a $z$-congruence. Also $\rho$ is a prime congruence, then $(f_1,g_1)\in \rho$ or $(f_2,g_2)\in \rho$. This implies either $A_1=E(f_1,g_1)\in E(\rho)$ or $A_2=E(f_2,g_2)\in E(\rho)$. Hence $E(\rho)$ is a prime $Z_\mathcal{A}$-filter on $X$.   
	\end{proof}
	
	The next theorem gives us a picture of the relation between prime congruences and $z$-congruences to some extent. 
	
	\begin{theorem}\label{P3, Th 4.15}
		For a $z$-congruence $\rho$ on $\mathcal{M}^{+}(X,\mathcal{A})$ the following are equivalent:
		
		\begin{itemize}
			\item[1.] $\rho$ is prime.
			
			\item[2.] $\rho$ contains a prime congruence.
			
			\item[3.] For all $f_1, f_2, g_1, g_2\in \mathcal{M}^{+}(X,\mathcal{A}), f_1f_2+g_1g_2=f_2g_2+f_2g_1$ implies that either $(f_1,g_1)\in \rho$ or $(f_2,g_2)\in \rho$.
			
			\item[4.] For all $f, g\in \mathcal{M}^{+}(X,\mathcal{A})$ there exists $A\in E(\rho)$ such that either $f\geq g$ or $g\geq f$ on $A$.
		\end{itemize}
	\end{theorem}
	
	\begin{proof}
		
		$(1)\Rightarrow (2):$ Trivial.
		
		$(2)\Rightarrow (3):$ Let $\rho$ contain a prime congruence $\sigma$ and $f_1, f_2, g_1, g_2\in \mathcal{M}^{+}(X,\mathcal{A}),$ such that $ f_1f_2+g_1g_2=f_2g_2+f_2g_1$. Then $(f_1f_2+g_1g_2, f_2g_2+f_2g_1)$ is a member of the diagonal congruence and hence a member of $\sigma$. Since $\sigma$ is prime, then either $(f_1, g_1)\in \rho_1$ or $(f_2, g_2)\in \sigma$. Since $\sigma\subseteq \rho$, so either $(f_1, g_1)\in \rho$ or $(f_2, g_2)\in \rho$.
		
		$(3)\Rightarrow (4):$ Let $f, g\in \mathcal{M}^{+}(X,\mathcal{A})$. We define $h_1=f-(f\wedge g)$ and $h_2=g-(g\wedge f)$. Then $h_1, h_2\in \mathcal{M}^{+}(X,\mathcal{A})$ and $h_1h_2=\boldsymbol{0}$. Thus $(h_1,\boldsymbol{0})\in \rho$ or $(\boldsymbol{0},h_2)\in \rho$ by $(3)$. This implies that $E(h_1,\boldsymbol{0})\in E(\rho)$ or $E(\boldsymbol{0}, h_2)\in E(\rho)$. Clearly, $f\geq g$ on $E(h_1,\boldsymbol{0})$ and $g\geq f$ on $E(h_2,\boldsymbol{0})$.
		
		$(4)\Rightarrow (1):$ Let $f_1, f_2, g_1, g_2\in \mathcal{M}^{+}(X,\mathcal{A})$ such that $(f_1f_2+g_1g_2, f_1g_2+f_2g_1)\in \rho$. Let $A=E(f_1f_2+g_1g_2, f_1g_2+f_2g_1)\in E(\rho)$. Set $h_1=\lvert f_1-g_1\rvert$ and $h_2=\lvert f_2-g_2\rvert$. Then $h_1, h_2\in \mathcal{M}^{+}(X,\mathcal{A}), E(f_1,g_1)=E(h_1,\boldsymbol{0})$ and $E(f_2,g_2)=E(h_2,\boldsymbol{0})$.
		
		By $(4)$ there exists $A_1\in E(\rho)$ such that $h_1\leq h_2$ or $h_2\leq h_1$ on $A_1$. We assume that $h_1\leq h_2$ on $A_1$. Then $h_1\leq h_2$ on $A\cap A_1$. Now $E(f_1f_2+g_1g_2, f_1g_2+f_2g_1)=E(f_1,g_1)\cup E(f_2,g_2)$. Therefore $A\cap A_1\subseteq E(h_1,\boldsymbol{0})=E(f_1,g_1)$. Then $E(f_1,g_1)\in E(\rho)$ as $A\cap A_1\in E(\rho)$ and $E(\rho)$ is an $\mathcal{A}$-filter on $X$. Thus $(f_1,g_1)\in \rho$ as $\rho$ is a $z$-congruence.
		
		Similarly $h_2\leq h_1$ on $A_1$ implies that $(f_2,g_2)\in \rho$. Hence, $\rho$ is a prime congruence.
	\end{proof}
	
	\begin{corollary}\label{prime-prime}
		If $\rho$ is a prime congruence on $\mathcal{M}^{+}(X,\mathcal{A})$, then $E(\rho)$ is a prime $\mathcal{A}$-filter on $X$.
	\end{corollary}
	
	\begin{proof}
		Let $\sigma= E^{-1}(E(\rho))$. Then $\sigma$ is a $z$-congruence and $\rho\subseteq \rho$. So $\sigma$ is prime by Theorem \ref{P3, Th 4.15}. Hence $E(\rho)=E(\sigma)$ is a prime $Z_\mathcal{A}$-filter on $X$ by Theorem \ref{P3, Th 4.15}.
	\end{proof}
	
	\begin{theorem}\label{prime cong}
		If $\mathfrak{F}$ is a prime $Z_\mathcal{A}$-filter on $X$, then $E^{-1}(\mathfrak{F})$ is a prime congruence on $\mathcal{M}^{+}(X,\mathcal{A})$.
	\end{theorem}
	
	\begin{proof}
		$E^{-1}(\mathfrak{F})$ is a congruence on $\mathcal{M}^{+}(X,\mathcal{A})$ by Proposition \ref{P3, prop 4.3}. Let $f_1, g_1, f_2, g_2 \in \mathcal{M}^{+}(X,\mathcal{A})$ such that $(f_1g_1+f_2g_2, f_1g_2+f_2g_1)\in E^{-1}(\mathfrak{F})$. Then $E(f_1g_1+f_2g_2, f_1g_2+f_2g_1)=E(f_1,g_1)\cup E(f_2,g_2) \in \mathfrak{F}$. Since $\mathfrak{F}$ is a prime $\mathcal{A}$-filter, therefore either $E(f_1,g_1)\in \mathfrak{F}$ or $E(f_2,g_2)\in \mathfrak{F}$. Thus either $(f_1,g_1)\in E^{-1}(\mathfrak{F})$ or $(f_2,g_2)\in E^{-1}(\mathfrak{F})$. Hence $E^{-1}(\mathfrak{F})$ is a prime congruence on $\mathcal{M}^{+}(X,\mathcal{A})$. 
	\end{proof}
	
	We have the following result as an easy consequence of Proposition \ref{P3, prop 4.3} and Theorem \ref{P3, Th4.5}.
	
	\begin{theorem}\label{max cong}
		The following statements hold for a measure space $(X,\mathcal{A})$.
		\begin{itemize}
			\item[a)] If $\rho$ is a maximal congruence on $\mathcal{M}^{+}(X,\mathcal{A})$ then $E(\rho)$ is an $Z_\mathcal{A}$-ultrafilter on $X$.
			
			\item[b)] If $\mathcal{U}$ is an $Z_\mathcal{A}$-ultrafilter on $X$ then $E^{-1}(\mathcal{U})$ is a maximal congruence on $\mathcal{M}^{+}(X,\mathcal{A})$.
		\end{itemize}
	\end{theorem}
	
	\begin{theorem}\label{max prime}
		Every maximal congruence on $\mathcal{M}^{+}(X,\mathcal{A})$ is a prime congruence.
	\end{theorem}
	
	\begin{proof}
		Let $\rho$ be a maximal congruence on $\mathcal{M}^{+}(X,\mathcal{A})$. Then $E(\rho)$ is an $\mathcal{A}$-ultrafilter on $X$. Every $\mathcal{A}$-ultrafilter is prime, therefore $E(\rho)$ is prime $\mathcal{A}$-filter on $X$. Since $\rho$ is a maximal congruence so it is a $z$-congruence. So we have $E^{-1}(E(\rho))=\rho$. Hence $\rho$ is prime by Theorem \ref{prime cong}.
	\end{proof}
	
	\begin{theorem}\label{prime maximal filter}
		Every prime $Z_\mathcal{A}$-filter on $X$ is an $Z_\mathcal{A}$-ultrafilter.
	\end{theorem}
	
	\begin{proof}
		Let $\mathfrak{F}$ be a prime $Z_\mathcal{A}$-filter on $X$. Then there exists an $\mathcal{A}$-ultrafilter $\mathcal{U}$ on $X$ such that $\mathfrak{F}\subseteq \mathcal{U}$. We show that $\mathfrak{F}=\mathcal{U}$. If $A\in \mathcal{U}\setminus \mathfrak{F}$, then $X=A\cup (X\setminus A)\in \mathfrak{F}$. This implies that $X\setminus A\in \mathfrak{F}$ as $A\notin \mathfrak{F}$ and $\mathfrak{F}$ is prime $Z_\mathcal{A}$-filter. We have $A, X\setminus A\in \mathcal{U}$. Thus $A\cap (X\setminus A)=\emptyset \in \mathcal{U}$, a contradiction. Hence $\mathfrak{F}=\mathcal{U}$ is an $Z_\mathcal{A}$-ultrafilter on $X$.
	\end{proof}
	
	\begin{theorem}\label{prime max}
		Let $\rho$ be a $z$-congruence on $\mathcal{M}^{+}(X,\mathcal{A})$. Then $\rho$ is maximal if and only if it is prime.
	\end{theorem}
	
	\begin{proof}
		If $\rho$ is maximal then it prime by Theorem \ref{max prime}.
		
		Conversely let, $\rho$ be prime. Then $E(\rho)$ is a prime $Z_\mathcal{A}$-filter by Theorem \ref{prime-prime}. By Theorem \ref{prime maximal filter}, $E(\rho)$ is an $\mathcal{A}$-ultrafilter. So $E^{-1}(E(\rho))$ is a maximal congruence on $\mathcal{M}^{+}(X,\mathcal{A})$ by Theorem \ref{max cong}. Since $\rho$ is a $z$-congruence, hence $E^{-1}(E(\rho))=\rho$ is maximal.
	\end{proof}

	We denote the intersection of all maximal congruences containing $(f,g)$ as $\mathfrak{M}(f,g)$. clearly, $\mathfrak{M}(f,g)$ is a $z$-congruence, for any $f,g\in \mathcal{M}^+(X,\mathcal{A})$
	Next we give an algebraic characterization of $z$-congruences on $\mathcal{M}^+(X,\mathcal{A})$. 
	
	\begin{proposition}
		A congruence $\rho$ on $\mathcal{M}^+(X,\mathcal{A})$ is a $z$-congruence if and only if $\mathfrak{M}(f,g)\subseteq \rho$ for every $(f,g)\in \rho$. 
	\end{proposition}
	\begin{proof}
		
		First, we observe that if $(h,k)$ belongs to every maximal congruence that $(f,g)$ belongs to, then $E(f,g)\subseteq E(h,k)$. Indeed, if we have $x\in E(f,g)$ and $x\notin E(h,k)$, then for this fixed point $x$, consider the fixed maximal ideal $M_x$ of $\mathcal{M}(X,\mathcal{A})$. Which forces  $E^{-1}(Z[M_x])$ to be a maximal congruence (cf. Theorem \ref{max cong}), where $Z[M_x]$ is a $Z_\mathcal{A}$-ultrafilter (cf. Theorem 2.7, \cite{A2020}). Then clearly $(f,g)\in E^{-1}(Z[M_x])$ but $(h,k)\notin E^{-1}(Z[M_x])$. Which is a contradiction. Therefore we have $E(f,g)\subseteq E(h,k)$.

		Now let $\rho$ be a $z$-congruence on $\mathcal{M}^+(X,\mathcal{A})$ and $(f,g)\in \rho$. For any $(h,k)\in \mathfrak{M}(f,g)$ we have $E(f,g)\subseteq E(h,k)$, which implies $(h,k)\in \rho$. Indeed, this follows from that fact that $E(\rho)$ is a $Z_\mathcal{A}$-filter and $\rho$ is a $z$-congruence. Therefore $\mathfrak{M}(f,g)\subseteq \rho$. 
		
		Conversely, it follows easily since each maximal congruence is a $z$-congruence. 
	\end{proof}
	
	\begin{corollary}
		Any $z$-congruence $\rho$ on $\mathcal{M}^+(X,\mathcal{A})$ is of the form $\rho=\underset{(f,g)\in \rho}{\bigvee_z}\mathfrak{M}(f,g)$.
	\end{corollary}

	\section{Structure space of $\mathcal{M}^{+}(X,\mathcal{A})$}\label{P3, Section 5}

	Let $Max(X,\mathcal{A})$ be the set of all maximal ideals of $\mathcal{M}(X,\mathcal{A})$ and for $f\in \mathcal{M}(X,\mathcal{A})$ set $\mathcal{M}_f=\{ M\in Max(X,\mathcal{A}): f\in M\}$. Then $\{\mathcal{M}_f: f\in \mathcal{M}(X,\mathcal{A})\}$ is base for closed sets for some topology on $Max(X,\mathcal{A})$. $Max(X,\mathcal{A})$ with this topology is called structure space of the ring $\mathcal{M}(X,\mathcal{A})$; See \cite{A2020}. 
	
    Let $S$ be a semiring and $\mathcal{M}Cong(S)$ be the set of all maximal congruences on $S$. For $a, b\in S$, set $\mathfrak{m}(a,b)=\{ \rho \in \mathcal{M}Cong(S): (a,b)\in \rho\}$. For $a, b, c, d\in S$, $\mathfrak{m}(a,b)\cup \mathfrak{m}(c,d)\subseteq \mathfrak{m}(ac+bd,ad+bc)$. If every maximal congruence on $S$ is prime, then the equality holds. Using this fact, we have the following theorem.
	
	\begin{theorem}[Theorem 2.9, \cite{A1993}]
		If each maximal congruence on $S$ is prime, then $\{ \mathfrak{m}(a,b): (a,b)\in S\times S\}$ is a base for closed sets of some topology on $\mathcal{M}Cong(S)$.
	\end{theorem}
	The set $\mathcal{M}Cong(S)$ with this topology is said to be the structure space of $S$, defined in \cite{S1990}.
	We now show that $Max(X,\mathcal{A})$ can be achieved via the positive cone of the ring $\mathcal{M}(X,\mathcal{A})$.
	
	\begin{theorem}\label{homeomorphism}
		$\mathcal{M}Cong(\mathcal{M}^{+}(X,\mathcal{A}))$ is homeomorphic to $Max(X,\mathcal{A})$.
	\end{theorem}
	
	\begin{proof}
		Let $\rho\in \mathcal{M}Cong(\mathcal{M}^{+}(X,\mathcal{A}))$. Then $E(\rho)$ is an $Z_\mathcal{A}$-ultrafilter on $X$ (cf. Theorem \ref{max cong}). Then $Z^{-1}[E(\rho)]$ is a maximal ideal in $\mathcal{M}(X,\mathcal{A})$ by Theorem 2.7 of \cite{A2020}. Define $\eta :\mathcal{M}Cong(\mathcal{M}^{+}(X,\mathcal{A}))\rightarrow Max(X,\mathcal{A})$ by $\eta (\rho)=Z^{-1}[E(\rho)]$ for every $\rho\in \mathcal{M}Cong(\mathcal{M}^{+}(X,\mathcal{A}))$.
		
		Let $\rho_1, \rho_2\in \mathcal{M}Cong(\mathcal{M}^{+}(X,\mathcal{A}))$ such that $\eta (\rho_1)=\eta (\rho_2)$. Then $Z^{-1}[E(\rho_1)]=Z^{-1}[E(\rho_2)]$. For any $Z_\mathcal{A}$-filter $\mathfrak{F}$, $ZZ^{-1}[\mathfrak{F}]=\mathfrak{F}$, so we have $E(\rho_1)=E(\rho_2)$. Also $\rho_1, \rho_2$ are $z$-congruences, therefore we have $\rho_1=\rho_2$. Hence $\eta$ is injective.
		
		Let $M$ be a maximal ideal of $\mathcal{M}(X,\mathcal{A})$. Then $Z[M]$ is an $Z_\mathcal{A}$-ultrafilter by Theorem 2.7 of \cite{A2020}. By Theorem \ref{max cong}, $E^{-1}Z[M]$ is a maximal congruence on $\mathcal{M}^{+}(X,\mathcal{A})$. Let $\rho=E^{-1}Z[M]$. Then $\eta(\rho)=Z^{-1}E(\rho)=M$. Hence $\eta$ is onto.
		
		Let $\mathfrak{m}(f,g)$ be a basic closed set in $\mathcal{M}Cong(\mathcal{M}^{+}(X,\mathcal{A}))$, where $f,g\in \mathcal{M}^{+}(X,\mathcal{A})$. Then $M\in \eta (\mathfrak{m}(f,g))\Leftrightarrow M=\eta (\rho), \rho\in \mathfrak{m}(f,g)\Leftrightarrow M=Z^{-1}E(\rho), (f,g)\in \rho\Leftrightarrow Z[M]=E(\rho), (f,g)\in \rho\Leftrightarrow E(f,g)=Z(f-g)\in E(\rho)=Z[M]\Leftrightarrow f-g\in M$ as $M$ is a $z$-ideal $\Leftrightarrow M\in \mathcal{M}_{f-g}$. Hence $\eta (\mathfrak{m}(f,g))=\mathcal{M}_{f-g}$ is a basic closed set in $Max(X,\mathcal{A})$.
		
		Let $\mathcal{M}_f$ be a basic closed set in $Max(X,\mathcal{A})$, where $f\in \mathcal{M}(X,\mathcal{A})$. Then $\rho \in \eta^{-1}(\mathcal{M}_f)\Leftrightarrow \eta (\rho)\in \mathcal{M}_f\Leftrightarrow f\in \eta (\rho)=Z^{-1}E(\rho)\Leftrightarrow Z(f)=Z(\lvert f\rvert)=E(\lvert f\rvert, \boldsymbol{0})\in E(\rho) \Leftrightarrow (\lvert f\rvert, \boldsymbol{0})\in \rho$ as $\rho$ is a $z$-congruence $\Leftrightarrow \rho \in m(\lvert f\rvert, \boldsymbol{0})$. Therefore $\eta^{-1}(\mathcal{M}_f)=m(\lvert f\rvert, \boldsymbol{0})$ is a basic closed set in $\mathcal{M}Cong(\mathcal{M}^{+}(X,\mathcal{A}))$. Hence $\eta$ is a homeomorphism.
	\end{proof}
	
	\begin{definition}
		A congruence $\rho$ is called \emph{fixed} if $\bigcap \{ E(f,g): (f,g)\in \rho\}$ is nonempty and \emph{free} otherwise. 
	\end{definition}
	
	
	\begin{theorem}\label{fixed}
		For a measurable space $(X,\mathcal{A})$, the set of all fixed maximal congruences of $\mathcal{M}^{+}(X,\mathcal{A})$ is the set $\{ \rho_x : x\in X\}$, where $\rho_x =\{ (f,g)\in \mathcal{M}^{+}(X,\mathcal{A})\times \mathcal{M}^{+}(X,\mathcal{A}): f(x)=g(x)\}$. Moreover, for two distinct points $x,y\in X$ we have $\rho_x \neq \rho_y$.
	\end{theorem}
	
	\begin{proof}
		For $x\in X$, we define $\phi_x: \mathcal{M}^{+}(X,\mathcal{A})\rightarrow \mathbb{R}_+$ by $\phi_x(f)=f(x)$ for all $f\in \mathcal{M}^{+}(X,\mathcal{A})$. Then clearly, $\phi_x$ is an onto semiring homomorphism with kernel $\phi_x$. Indeed, $\phi_x=\{ (f,g):\mathcal{M}^{+}(X,\mathcal{A})\times \mathcal{M}^{+}(X,\mathcal{A}): \phi_x(f)=\phi_x(g)\}=\{ (f,g):\mathcal{M}^{+}(X,\mathcal{A})\times \mathcal{M}^{+}(X,\mathcal{A}): f(x)=g(x)\}=\rho_x$. Therefore $\mathcal{M}^{+}(X,\mathcal{A})/\rho_x$ is isomorphic to $\mathbb{R}_+$. Therefore $\mathcal{M}^{+}(X,\mathcal{A})/\rho_x$ is congruence free as $\mathbb{R}_+$ is so. Hence $\rho_x$ is a maximal congruence on $\mathcal{M}^{+}(X,\mathcal{A})$ and also clearly it is fixed.
		
		Conversely let, $\rho$ be fixed maximal congruence of $\mathcal{M}^{+}(X,\mathcal{A})$ and $x\in \bigcap \{ E(f,g): (f,g)\in \rho\}$. Then $\rho\subseteq \rho_x$ and hence by maximality of $\rho$ we have $\rho=\rho_x$.
		
		Since $\mathcal{A}$ separate points of $X$. Then for $x\neq y$ in $X$ there exists $A\in \mathcal{A}$ such that $x\in A$ and $y\notin A$. Then $A= E(f,\boldsymbol{0})$, where $f=\chi_{A^c}\in \mathcal{M}^{+}(X,\mathcal{A})$. So we have $(f,\boldsymbol{0})\in \rho_x$ but $(f,\boldsymbol{0})\notin \rho_y$. Hence $\rho_x\neq \rho_y$.
	\end{proof}
	
Next, we try to characterize the measurable space $(X,\mathcal{A})$ for which each maximal congruence on $\mathcal{M}^+(X,\mathcal{A})$ is fixed. We need the following definitions. 

\begin{definition}(cf. Definition 7.15 of \cite{Davey})
    Let $L$ be a complete lattice, and let $a$ be an element of $L$. Then $a$ is called \emph{compact} if $a\leq \bigvee X$, for some $X\subseteq L$, implies that $a\leq \bigvee S$ for some finite $S\subseteq X$. In particular, if the top element of $L$ is compact, then we call $L$ a \emph{compact lattice}. 
\end{definition}

 \begin{definition}(\cite{Estaji})
     A measurable space $(X,\mathcal{A})$ is said to be a \emph{compact measurable space} if $\mathcal{A}$ is a compact lattice. 
 \end{definition}
	
	\begin{theorem}\label{fixed cong}
		Let $(X,\mathcal{A})$ be a measurable space. Then the following are equivalent.
		
		\begin{itemize}
			\item[1)] Each maximal ideal of $\mathcal{M}(X,\mathcal{A})$ is fixed.
			
			\item[2)] Each maximal congruence on $\mathcal{M}^{+}(X,\mathcal{A})$ is fixed.
			
			\item[3)] $\mathcal{A}$ is a finite $\sigma$-algebra on $X$.

                \item[4)] $(X,\mathcal{A})$ is a compact measurable space. 

                \item[5)] $\mathcal{M}(X,\mathcal{A})=\mathcal{M}^*(X,\mathcal{A})=\{f\in \mathcal{M}(X,\mathcal{A})\colon f\ \text{is bounded on}\ X\}$.  
		\end{itemize}
	\end{theorem}
	
	\begin{proof}
		$(1)\Leftrightarrow (3)\Leftrightarrow (5):$ Follows from Theorem 2.13 of \cite{A2020}.

            $(1)\Leftrightarrow (4)$ Follows from Proposition 4.11 of \cite{Estaji}.

		$(1)\Rightarrow (2):$ Let $\rho$ be a maximal congruence on $\mathcal{M}^{+}(X,\mathcal{A})$. Then $E(\rho)$ is an $Z_{\mathcal{A}}$-ultrafiler on $X$. So $Z^{-1}[E(\rho)]$ is a maximal ideal of $\mathcal{M}(X,\mathcal{A})$. Thus $Z^{-1}[E(\rho)]=\{ h\in \mathcal{M}(X,\mathcal{A}): h(x)=0\}$ for some $x\in X$. Let $(f,g)\in \rho$. Then $E(f,g)=Z(f-g)\in E(\rho)$. Therefore $f-g\in Z^{-1}[E(\rho)]$. So $(f-g)(x)=0$ i.e., $x\in E(f,g)$. Hence $\rho$ is a fixed congruence on $\mathcal{M}^{+}(X,\mathcal{A})$.
		
		$(2)\Rightarrow (1):$ Let $M$ be any maximal ideal of $\mathcal{M}(X,\mathcal{A})$. Then $Z[M]$ is an $Z_\mathcal{A}$-utrafilter on $X$. Thus $E^{-1}(Z[M])$ is a maximal congruence on $\mathcal{M}^{+}(X,\mathcal{A})$. By $(2)$ $E^{-1}(Z[M])$ is fixed. So by Theorem \ref{fixed} $E^{-1}(Z[M])=\rho_x$ for some $x\in X$. Let $f\in M$. Then $Z(f)=E(\lvert f\rvert,\boldsymbol{0})\in Z[M]$. Thus $(\lvert f\rvert,\boldsymbol{0})\in E^{-1}(Z[M])=\rho_x$. Therefore $f(x)=\lvert f\rvert(x)=\boldsymbol{0}(x)=0$. Hence $M$ is a fixed maximal ideal of $\mathcal{M}(X,\mathcal{A})$.

            $(3)\Rightarrow (4)$ Trivially true. 

            $(4)\Rightarrow (3)$ Suppose $\mathcal{A}$ is an infinite $\sigma$-algebra on $X$. For each $x\in X$, we define $x^{\curlyvee}=\bigcap \{B\in \mathcal{A}\colon x\in B\}$. Since $\mathcal{A}$ is a complete lattice, we have $x^{\curlyvee}\in \mathcal{A}$ for all $x\in X$. It is easy to see that $x^{\curlyvee}\cap y^{\curlyvee}=\emptyset$ or $x^{\curlyvee}=y^{\curlyvee}$, for two distinct points $x,y\in X$. Let $\mathfrak{X}=\{x^{\curlyvee}\colon x\in X\}$.
            It is clear that any $B\in \mathcal{A}$ can be written as unions of elements of $\mathfrak{X}$. Since $\mathcal{A}$ is infinite, then $\mathfrak{X}$ has infinite cardinality. Now let, $\{x_1^{\curlyvee}, x_2^{\curlyvee},\cdots\}$ be a countable infinite subset of $\mathfrak{X}$. Let $A=\bigcup_i x_i^{\curlyvee}$. Then $A\in \mathcal{A}$. By our hypothesis, $(X, \mathcal{A})$ is a compact measurable space. Then there exists a finite subcollection $\{{x_i}_j^{\curlyvee}\}_j$ of $\{x_i^{\curlyvee}\}_i$ such that $A=\bigcup_j {x_i}_j^{\curlyvee}$. This contradicts the fact that each member of $\{x_i^{\curlyvee}\}_i$ is disjoint.

	\end{proof}

	\begin{remark}
	    From the definition of a $\sigma$-algebra, it is clear that a $\sigma$-algebra is a $\sigma$-frame (A $\sigma$-frame $L$ is a lattice with countable joins $\bigvee_n$, finite meets $\wedge$, a top element $\top$ and a bottom element $\perp$ such that $x\wedge \bigvee_n x_n=\bigvee_n (x\wedge x_n)$, for $n\in \mathbb{N}$, for all $x, x_n\in L$). A $\sigma$-frame $S$ is said to be \emph{pseudocompact} if every $\sigma$-frame maps $\phi\colon \mathfrak{L}(\mathbb{R})\rightarrow S$ is a bounded map (cf. Definition 3, \cite{Gilmour}). Here $\mathcal{L}(\mathbb{R})$ is the $\sigma$-frame (frame) of real numbers, which is isomorphic to $\mathcal{O}(\mathbb{R})$, the lattice of open sets of $\mathbb{R}$. Therefore a $\sigma$-algebra $\mathcal{A}$ is pseudocompact if every $f\colon \mathcal{O}(\mathbb{R})\rightarrow \mathcal{A}$ is bounded. Unlike the classical case of rings of continuous functions, that is, a space $X$ is pseudocompact if and only if $C(X)=C^*(X)$, Theorem \ref{fixed cong} is unable to capture pseudocompactness of $(X,\mathcal{A})$. Now, we argue why this phenomenon occurs. We show that any $\sigma$-algebra is a regular $\sigma$-frame (for definition, see the Background section of \cite{Gilmour}). Consider any nonempty element $B$ of $\mathcal{A}$. If $B$ is an \textit{atom}, then vacuously $B$ is a regular element of the $\sigma$-frame $\mathcal{A}$. Now let $B$ be a non-atom of $\mathcal{A}$. Then there exists $\{E_n\colon n\in \mathbb{N}\}$,  a pairwise disjoint countable nonempty members of $\mathcal{A}$ defined as in Lemma 2.12 of \cite{A2020}. It is easy to see that $E_{i+1}\prec E_i$ and $B=\bigvee_i E_i$, $i=1,2,\cdots$. Therefore every non-atom is also a regular element. We conclude that $\mathcal{A}$ is a regular $\sigma$-frame. Now under the \textit{Axiom of countable dependant choice}, any regular $\sigma$-frame is completely regular. Then by Corollary 2 of \cite{Gilmour}, compactness and pseudocompactness coincide in $\mathcal{A}$. 
	\end{remark}
	
	\begin{proposition}\label{P3, Prop5.9}
	    For every compact measurable space $(X,\mathcal{A})$, 
     \begin{center}
         $(\mathcal{M}Cong(\mathcal{M}^+(X,\mathcal{A}), \{\mathfrak{m}(f,g)\colon f,g\in \mathcal{M}^+(X,\mathcal{A})\})$  is $T$-measurable. 
     \end{center} 
	\end{proposition}
	
	\begin{proof}
	    We know from Theorem \ref{fixed cong} that each maximal congruence on $\mathcal{M}^+(X,\mathcal{A})$ is of the form $\rho_x$ for each $x\in X$ (cf. Theorem \ref{fixed}). It follows that $\mathfrak{m}(f,g)=\{\rho_x\colon (f,g)\in \rho_x\}=\{\rho_x\colon x\in E(f,g)\}$. Therefore $\mathfrak{m}(f,g)^c=\{\rho_x\colon x\in E^c(f,g)\}=\{\rho_x\colon x\in E(\chi_{E(f,g)},\boldsymbol{0})\}=\mathfrak{m}(\chi_{E(f,g)},\boldsymbol{0})$. Moreover $\bigcup_{n\in \mathbb{N}} \mathfrak{m}(f_n,g_n)=\bigcup_{n\in \mathbb{N}} \{\rho_x\colon x\in E(f_n,g_n)\}=\{\rho_x\colon x\in \bigcup_{n\in \mathbb{N}} E(f_n,g_n)\}=\mathfrak{m}(\chi_{\cap_{n\in \mathbb{N}}E^c(f_n,g_n)}, \boldsymbol{0})$. Therefore $(\mathcal{M}Cong(\mathcal{M}^+(X,\mathcal{A}), \{\mathfrak{m}(f,g)\colon f,g\in \mathcal{M}^+(X,\mathcal{A})\})$ is a measurable space and by Theorem \ref{fixed} it is a $T$-measurable space. 
	\end{proof}

 \begin{definition}(\cite{Estaji})
     Let $(X,\mathcal{A})$ and $(Y,\mathcal{B})$ be two measurable spaces. We say that $(X,\mathcal{A})$ and $(Y,\mathcal{B})$ are homeomorphic if there exists a one-one and onto function $f\colon X\rightarrow Y$ such that $A\in \mathcal{A}$ if and only if $f(A)\in \mathcal{B}$, for every $A\in \mathcal{A}$.
 \end{definition}
	When $(X,\mathcal{A})$ is homeomorphic to $(Y,\mathcal{B})$ we will simply write $X\cong Y$.
 \begin{theorem}\label{P3, Th5.11}
     For every compact measurable space $(X,\mathcal{A})$, 
     
     $X\cong \mathcal{M}Cong(\mathcal{M}^+(X,\mathcal{A}))$ as measurable spaces.
 \end{theorem}
	\begin{proof}
	    We define $\phi\colon X\rightarrow \mathcal{M}Cong(\mathcal{M}^+(X,\mathcal{A}))$ by $\phi(x)=\rho_x$. From Theorem \ref{fixed cong} it is clear that $\phi$ is one-one and onto, Moreover $\phi(E(f,g))=\{\rho_x\colon x\in E(f,g)\}=\mathfrak{m}(f,g)$ and $\phi^{-1}(\mathfrak{m}(f,g))=E(f,g)$. Therefore, $\phi$ is a homeomorphism.
	\end{proof}
 The following corollary is a direct consequence of Theorem \ref{homeomorphism} and Theorem \ref{P3, Th5.11}.
 \begin{corollary}\label{P3, Cor5.12}
     If $(X,\mathcal{A})$ and $(Y,\mathcal{B})$ are two compact measurable spaces, then $X\cong Y$ as measurable spaces if and only if $\mathcal{M}^+(X, \mathcal{A})\cong \mathcal{M}^+(Y,\mathcal{B})$ as semirings. 
 \end{corollary}

	\section{Real maximal congruences on $\mathcal{M}^{+}(X,\mathcal{A})$}\label{P3, Section 6}
	
	Here we initiate a study of quotients of $\mathcal{M}^+(X,\mathcal{A})$ via some important class of congruences on it. Our goal is to give an alternative description of realcompact measurable spaces, in view of $\mathcal{M}^+(X,\mathcal{A})$. The following lemma infers that the class of $z$-congruences is ideal to consider while sculpting quotients of $\mathcal{M}^+(X,\mathcal{A})$.
	
	\begin{lemma}\label{P3, Lem6.1}
		Every $z$-congruence on $\mathcal{M}^{+}(X,\mathcal{A})$ is convex.
	\end{lemma}
	
	\begin{proof}
		Let $\rho$ be a $z$-congruence on $\mathcal{M}^{+}(X,\mathcal{A})$. Let $f,g,f_1,g_1\in \mathcal{M}^{+}(X,\mathcal{A})$ such that $f\leq f_1\leq g_1\leq g$ and $(f,g)\in \rho$. Then $E(f,g)\subseteq E(f_1,g_1)$, $E(f,g)\in E(\rho)$ and $E(\rho)$ is an $Z_\mathcal{A}$-filter $\Rightarrow$ $E(f_1,g_1)\in E(\rho)$. Also $E(\rho)$ is a $z$-congruence. Therefore $(f_1,g_1)\in \rho$. Hence $\rho$ is a convex congruence. 
	\end{proof}
	
	\begin{theorem}\label{ordering}
		Let $\rho$ be a $z$-congruence on $\mathcal{M}^{+}(X,\mathcal{A})$ and $f,g\in \mathcal{M}^{+}(X,\mathcal{A})$. Then $\rho (f)\leq \rho (g)$ if and only if $f\leq g$ on some member of $E(\rho)$.
	\end{theorem}
	
	\begin{proof}
		Let $f, g\in \mathcal{M}^{+}(X,\mathcal{A})$ such that $\rho (f)\leq \rho (g)$. Then by Theorem \ref{P3, Th2.4} and Lemma \ref{P3, Lem6.1}, there exists $(h,k)\in \rho$ such that $f+h\leq g+k$. We have $f\leq g$ on $E(h,k)\in E(\rho)$.
		
		Conversely, let $f\leq g$ on some $A\in E(\rho)$. Then $A=E(f_1,g_1)$ for some $(f_1,g_1)\in \rho$. Set $h=(f-g)\vee \boldsymbol{0}$. Then $h\in \mathcal{M}^{+}(X,\mathcal{A})$ and $A\subseteq E(h,\boldsymbol{0})$. Since $E(\rho)$ is a $Z_{\mathcal{A}}$-filter, therefore $E(h,\boldsymbol{0})\in E(\rho)$. This implies that $(h,\boldsymbol{0})\in \rho$ as $\rho$ is a $z$-congruence. Also $f+\boldsymbol{0}\leq g+h$. Hence $\rho (f)\leq \rho (g)$.
	\end{proof}
	
	\begin{theorem}
		Let $f,g\in \mathcal{M}^{+}(X,\mathcal{A})$ and $\rho$ be a maximal congruence on $\mathcal{M}^{+}(X,\mathcal{A})$. Then $\rho (f)<\rho (g)$ if and only if there exists $Z\in E(\rho)$ such that $f<g$ on $Z$. 
	\end{theorem}
	
	\begin{proof}
		Let $f<g$ on some $Z\in E(\rho)$. Then $E(f,g)\cap Z=\emptyset$. Therefore $E(f,g)\notin E(\rho)$. This implies that $(f,g)\notin \rho$. So, $\rho(f)\neq \rho (g)$. By the above theorem, we get $\rho (f)<\rho (g)$.
		
		Conversely let, $\rho (f)<\rho (g)$. Since every maximal congruence is a $z$-congruence, then by Theorem \ref{ordering} we have $f\leq g$ on some $Z_1\in E(\rho)$. Again $\rho (f)\neq \rho (g)$, this implies $(f,g)\notin \rho$.  Therefore $E(f,g)\notin E(\rho)$ as $\rho$ is a $z$-congruence. As $E(\rho)$ is an $Z_{\mathcal{A}}$-ultrafiler then there exists $Z_2\in E(\rho)$ such that $E(f,g)\cap Z_2=\emptyset$. Then $f<g$ on $Z_1\cap Z_2\in E(\rho)$. This completes the proof.
	\end{proof}
	
	\begin{theorem}
		If $\rho$ is a maximal congruence on $\mathcal{M}^{+}(X,\mathcal{A})$, then $\mathcal{M}^{+}(X,\mathcal{A})/\rho$ is totally ordered semiring.
	\end{theorem}
	
	\begin{proof}
		Let $f,g \in \mathcal{M}^{+}(X,\mathcal{A})$. Set $A_1=\{ x\in X: f(x)\leq g(x)\}$, $A_2=\{ x\in X: g(x)\leq f(x)\}$. Then $A_1, A_2\in \mathcal{A}$ and $A_1\cup A_2=X$. Since $\rho$ is a maximal congruence, it is prime. Therefore $E(\rho)$ is a prime $Z_{\mathcal{A}}$-filter. Thus either $A_1\in E(\rho)$ or $A_2\in E(\rho)$. Now $f\leq g$ on $A_1$ and $g\leq f$ on $A_2$. Then by Theorem \ref{ordering} either $\rho (f)\leq \rho (g)$ or $\rho (g)\leq \rho (f)$. Hence $\mathcal{M}^{+}(X,\mathcal{A})/\rho$ is totally ordered semiring.
	\end{proof}
	
	We can easily prove the following result.
	\begin{theorem}\label{orderpreserving}
		Let $\rho$ be maximal congruence on $\mathcal{M}^{+}(X,\mathcal{A})$. then the mapping $\phi : \mathbb{R}_+\rightarrow \mathcal{M}^{+}(X,\mathcal{A})/\rho$, defined by $\phi (r)=\rho (\boldsymbol{r})$ is an ordered preserving isomorphism from $\mathbb{R}_+$ into $\mathcal{M}^{+}(X,\mathcal{A})/\rho$.
	\end{theorem}
	
	\begin{definition}
		Let $\rho$ be a maximal congruence on $\mathcal{M}^{+}(X,\mathcal{A})$. Then $\rho$ is called real if $\phi$ is onto where $\phi$ is defined in Theorem \ref{orderpreserving}. A maximal congruence is said to be hyper-real if it is not real.
	\end{definition}
	
	\begin{theorem}\label{real}
		A maximal congruence $\rho$ on $\mathcal{M}^{+}(X,\mathcal{A})$ is real if and only if $\mathcal{M}^{+}(X,\mathcal{A})/\rho$ is isomorphic to $\mathbb{R}_+$.
	\end{theorem}
	
	\begin{proof}
		If $\rho$ is real then by the definition of real congruence $\mathcal{M}^{+}(X,\mathcal{A})/\rho$ is isomorphic to $\mathbb{R}_+$.
		
		Conversely let, $\mathcal{M}^{+}(X,\mathcal{A})/\rho$ be isomorphic to $\mathbb{R}_+$ and $\psi$ is an isomorphism form $\mathcal{M}^{+}(X,\mathcal{A})/\rho$ onto $\mathbb{R}_+$. Then $\phi\circ \psi $ is an isomorphism from $\mathbb{R}_+$ into $\mathbb{R}_+$. But only non-zero isomorphism from $\mathbb{R}_+$ into $\mathbb{R}_+$ is the identity map. Therefore $\phi \circ \psi$ is the identity map. So $\phi$ is onto. Hence $\rho$ is a real maximal congruence. 
	\end{proof}
	
	\begin{theorem}\label{P3, Th6.8}
		For each $x\in X$, the fixed congruence $\rho_x=\{ (f,g)\in \mathcal{M}^{+}(X,\mathcal{A})\times \mathcal{M}^{+}(X,\mathcal{A}): f(x)=g(x)\}$ on $\mathcal{M}^{+}(X,\mathcal{A})$ is real.
	\end{theorem}
	
	\begin{proof}
		Follows from Theorem \ref{fixed} and Theorem \ref{real}.
	\end{proof}
	
	Lemma \ref{inverse}, Theorem \ref{cofinal}, Theorem \ref{infinte}, \ref{closed under countable intersection} follows arguing similarly as in the proof of Lemma 4.9, Theorem 4.8, Theorem 4.10,
 and Theorem 4.11 respectively in \cite{A1995}.
	
	\begin{lemma}\label{inverse}
		For any maximal congruence $\rho$ on $\mathcal{M}^{+}(X,\mathcal{A})$ each non-zero element in $\mathcal{M}^{+}(X,\mathcal{A})/\rho$ has multiplicative inverse.
	\end{lemma}
	
	\begin{theorem}\label{cofinal}
		A maximal congruence $\rho$ on $\mathcal{M}^{+}(X,\mathcal{A})$ is real if and only if the set $\{ \rho (\boldsymbol{n}): n\in \mathbb{N}\}$ is cofinal in $\mathcal{M}^{+}(X,\mathcal{A})/\rho$.
	\end{theorem}
	
	An element $a$ in a totally ordered semiring $S$ is called \emph{infinitely large} if $a\geq n$ for all $n\in \mathbb{N}$.
	
	\begin{theorem}\label{infinte}
		Let $\rho$ be a maximal congruence on $\mathcal{M}^{+}(X,\mathcal{A})$ and $f\in \mathcal{M}^{+}(X,\mathcal{A})$. Then the following statements are equivalent:
		
		\begin{itemize}
			\item[1)] $\rho (f)$ is infinitely large.
			
			\item[2)] For all $n\in \mathbb{N}$, $Z_n=\{ x\in X: f(x)\geq n\} \in E(\rho)$.
			
			\item[3)] For all $n\in \mathbb{N}, (f\wedge \boldsymbol{n},\boldsymbol{n})\in \rho$.
			
			\item[4)] $f$ is unbounded on each member of $E(\rho)$.
		\end{itemize} 
	\end{theorem}
	
	\begin{theorem}\label{closed under countable intersection}
		A maximal congruence $\rho$ is real if and only if $E(\rho)$ is closed under countable intersection.
	\end{theorem}
	
	\begin{theorem}\label{real congruence}
		A maximal congruence $\rho$ on $\mathcal{M}^{+}(X,\mathcal{A})$ is real if and only if $Z^{-1}[E(\rho)]$ is a real maximal ideal in $\mathcal{M}(X,\mathcal{A})$.
	\end{theorem}
	
	\begin{proof}
		A maximal ideal $M$ in $\mathcal{M}(X,\mathcal{A})$ is real if and only if $Z[M]$ is closed under countable intersection by Theorem 2.2\cite{AA2020}. Therefore the maximal ideal $Z^{-1}[E(\rho)]$ is real if and only if $ZZ^{-1}E(\rho)=E(\rho)$ is closed under countable intersection if and only if $\rho $ is real by Theorem \ref{closed under countable intersection}. This completes the proof.
	\end{proof}
	
	Let $\mathcal{RM}Cong(\mathcal{M}^{+}(X,\mathcal{A}))=\{ \rho\in \mathcal{M}Cong(\mathcal{M}^{+}(X,\mathcal{A})): \rho$ is a real maximal congruence$\}$ and $RMax(X,\mathcal{A})$ be the set of all real maximal ideals  of $\mathcal{M}(X,\mathcal{A}))$. 
 Topologically, $\mathcal{RM}Cong(\mathcal{M}^{+}(X,\mathcal{A}))$ is a subspace of 
 
 $\mathcal{M}Cong(\mathcal{M}^{+}(X,\mathcal{A}))$ with basic closed sets 
 
 $\mathfrak{m}^R(f,g)=\{\rho\in \mathcal{RM}Cong(\mathcal{M}^{+}(X,\mathcal{A}))\colon (f,g)\in \rho\}$.
	
	\begin{theorem}
		The map $\Tilde{\eta}: \mathcal{RM}Cong(\mathcal{M}^{+}(X,\mathcal{A}))\rightarrow RMax(X,\mathcal{A})$ is a homeomorphism, where $\Tilde{\eta}$ is the restriction map $\eta\lvert_{\mathcal{RM}Cong(\mathcal{M}^{+}(X,\mathcal{A}))}$ (cf. Theorem \ref{homeomorphism}).
	\end{theorem}
	
	\begin{proof}
		For $\rho \in \mathcal{RM}Cong(\mathcal{M}^{+}(X,\mathcal{A})), \Tilde{\eta}(\rho)\in RMax(X,\mathcal{A})$ by Theorem \ref{real congruence}.
		In view of Theorem \ref{homeomorphism} it is only to show that $\Tilde{\eta}$ is onto. Let $M\in RMax(X,\mathcal{A})$. Then $Z[M]$ is an $Z_{\mathcal{A}}$-utrafilter. Then $Z[M]=E(\rho)$ for a unique maximal congruence $\rho$ on $\mathcal{M}^{+}(X,\mathcal{A})$. Since $M$ is real, $Z[M]=E(\rho)$ is closed under countable intersection. Then by Theorem \ref{closed under countable intersection} $\rho$ is real i.e., $\rho\in \mathcal{RM}Cong(\mathcal{M}^{+}(X,\mathcal{A}))$. Also $\Tilde{\eta}(\rho)=Z^{-1}E(\rho)=Z^{-1}Z[M]=M$. Thus $\Tilde{\eta}$ is onto. Moreover if $\mathcal{M}_f^R$ is a basic closed set in $RM(X,\mathcal{A})$, then $\Tilde{\eta}(\mathfrak{m}^R(f,g))=\mathcal{M}_{f-g}^R$ for $f,g\in \mathcal{M}^+(X,\mathcal{A})$ and $\Tilde{\eta}^{-1}(\mathcal{M}_h^R)=\mathfrak{m}^R(|h|,\boldsymbol{0})$ for $h\in \mathcal{M}(X,\mathcal{A})$. This completes the proof.
	\end{proof}

        The collection $\mathcal{RM}Cong(\mathcal{M}^{+}(X,\mathcal{A}))$ can be made into a measurable space also.  
        \begin{theorem}
            For a measurable space $(X,\mathcal{A})$, 
            
            $(\mathcal{RM}Cong(\mathcal{M}^+(X,\mathcal{A})), \{\mathfrak{m}^R(f,g)\colon f,g\in \mathcal{M}^+(X,\mathcal{A})\})$ forms a $T$-measurable space.  
        \end{theorem}
	\begin{proof}
	    To show that $\{\mathfrak{m}^R(f,g)\colon f,g\in \mathcal{M}^+(X,\mathcal{A})\}$ is a $\sigma$-algebra on $\mathcal{RM}Cong(\mathcal{M}^+(X,\mathcal{A})$. First we observe that $\mathfrak{m}^R(\boldsymbol{1},\boldsymbol{0})=\emptyset$. Indeed, no proper congruence contains the identity pair $(\boldsymbol{1},\boldsymbol{0})$. Next, for any $\mathfrak{m}^R(f,g)$, the complement is 
     
     $\mathfrak{m}^R(f,g)^c=\mathcal{RM}Cong(\mathcal{M}^+(X,\mathcal{A}))\setminus \{\mathfrak{m}^R(f,g)\}=\mathfrak{m}^R(\chi_{Z(f)\cup Z(g)},\boldsymbol{0})$. Indeed the twisted product $(f,g)\cdot_t (\chi_{Z(f)\cup Z(g)},\boldsymbol{0})\in \rho$, for all $\rho\in \mathfrak{m}^R(f,g)^c$ and $\rho$ does not contain $(f,g)$. Since every maximal congruence is prime, $(\chi_{Z(f)\cup Z(g)},\boldsymbol{0})\in \rho$. Hence $\mathfrak{m}^R(f,g)^c\subseteq \mathfrak{m}^R(\chi_{Z(f)\cup Z(g)},\boldsymbol{0})$ and reverse inclusion follows easily. Finally, let $\rho\in \bigcup_{n\in \mathbb{N}}\mathfrak{m}^R(f_n,g_n)$. Then $\rho\in \mathfrak{m}^R(f_n,g_n)$ for some $n\in \mathbb{N}$, so $(f_n,g_n)\in \rho$. Therefore $E(\chi_{\bigcap_{k\in \mathbb{N}}E^c(f_k,g_k)},\boldsymbol{0})\supseteq E(f_n,g_n)$. Hence $\bigcup_{n\in \mathbb{N}}\mathfrak{m}^R(f_n,g_n)\subseteq \mathfrak{m}^R(\chi_{\bigcap_{k\in \mathbb{N}}E^c(f_k,g_k)},\boldsymbol{0})$. 
     The reverse inclusion follows easily. Therefore $\{\mathfrak{m}^R(f,g)\colon f,g\in \mathcal{M}^+(X,\mathcal{A})\}$ is a $\sigma$-algebra on $\mathcal{RM}Cong(\mathcal{M}^+(X,\mathcal{A})$. Moreover, since all fixed maximal congruences are real (cf. Theorem \ref{P3, Th6.8}), the $\sigma$-algebra $\{\mathfrak{m}^R(f,g)\colon f,g\in \mathcal{M}^+(X,\mathcal{A})\}$ separates points (cf. Theorem \ref{fixed}). 
	\end{proof}
\begin{lemma}
     There is a one-one correspondence between real maximal congruences on $\mathcal{M}^+(X,\mathcal{A})$ and real maximal ideals of $\mathcal{M}(X,\mathcal{A})$.
 \end{lemma}
 \begin{proof}
     It is evident from Proposition \ref{P3, Prop4.12} that every  maximal congruence on $\mathcal{M}^+(X,\mathcal{A})$ is of the form $k^\nabla$, where $k$ is a maximal congruence on $\mathcal{M}(X,\mathcal{A})$. Now by Theorem 7.1 of \cite{Hebisch}, the quotient ring $\mathcal{M}(X,\mathcal{A})/k$ is isomorphic to the ring of differences of the quotient semiring $\mathcal{M}^+(X,\mathcal{A})/k^{\nabla}$. Therefore $k$ is a real maximal congruence on $\mathcal{M}(X, A)$ if and only if $k^\nabla$ is a real maximal congruence on $\mathcal{M}^+(X,\mathcal{A})$. Since $\mathcal{M}(X,\mathcal{A})/k\cong \mathcal{M}(X,\mathcal{A})/I_k$, this completes the proof. 
 \end{proof}
 
 In view of the preceding proposition, we can define realcompactness as follows: A measurable space $(X,\mathcal{A})$ is said to be realcompact if every real maximal congruence on $\mathcal{M}^+(X,\mathcal{A})$ is fixed. 

 Let us denote the set of all fixed maximal congruences as $\mathcal{FM}Cong(\mathcal{M}^+(X,\mathcal{A})$. Then by Theorem \ref{P3, Th6.8}, $\mathcal{FM}Cong(\mathcal{M}^+(X,\mathcal{A}))\subseteq \mathcal{RM}Cong(\mathcal{M}^+(X,\mathcal{A}))$. Equality holds for realcompact spaces. 
\begin{lemma}\label{P3, lem6.17}
    For any measurable space $(X,\mathcal{A})$, $\mathcal{FM}Cong(\mathcal{M}^+(X,\mathcal{A})\cong X$.
\end{lemma}
\begin{proof}
    Follows from Proposition \ref{P3, Prop5.9} and Theorem \ref{P3, Th5.11}.  
\end{proof}
 \begin{theorem}\label{P3, Th6.18}
     For a realcompact space $(X,\mathcal{A})$ the measurable spaces $(RMax(X,\mathcal{A}), \{\mathcal{B}_f^R\colon f\in f\in \mathcal{M}(X,\mathcal{A})\})$ (cf. Theorem 2.4, \cite{AA2020}) and 
     
     $(\mathcal{RM}Cong(\mathcal{M}^+(X,\mathcal{A})),\{\mathfrak{m}^R{(f,g)}\colon f,g\in \mathcal{M}^+(X,\mathcal{A})\})$ are homeomorphic. 
 \end{theorem}
 \begin{proof}
     Since $(X,\mathcal{A})$ is a realcompact space, $\mathcal{RM}Cong(\mathcal{M}^+(X,\mathcal{A}))=\mathcal{FM}Cong(\mathcal{M}^+(X,\mathcal{A}))$ and $RMax(X,\mathcal{A})=FMax(X,\mathcal{A})$, where $FMax(X,\mathcal{A})$ is the set of all fixed maximal ideals of $\mathcal{M}(X,\mathcal{A})$. By Lemma \ref{P3, lem6.17} and Theorem 2.6 of \cite{AA2020}, we conclude that $\mathcal{RM}Cong(\mathcal{M}^+(X,\mathcal{A}))\cong X\cong RMax(X,\mathcal{A})$.
 \end{proof}
 \begin{theorem}\label{P3, Th6.19}
     Two realcompact measurable spaces $(X,\mathcal{A})$ and $(Y,\mathcal{B})$ are homeomorphic if and only if $\mathcal{M}^+(X,\mathcal{A})$ and $\mathcal{M}^+(Y,\mathcal{B})$ are isomorphic as semirings.  
 \end{theorem}
 \begin{proof}
     If $(X,\mathcal{A})$ and $(Y,\mathcal{B})$ are homeomorphic, then $\mathcal{M}(X,\mathcal{A})$ and $\mathcal{M}(Y,\mathcal{B})$ are isomorphic as semirings. Since every ring homomorphism $\phi$ between $\mathcal{M}(X,\mathcal{A})$ and $\mathcal{M}(Y,\mathcal{B})$ is also a lattice homomorphism, the homomorphism preserves order. Indeed, if $f\leq g$ in $\mathcal{M}(X,\mathcal{A})$, then $g-f$ is a square and so $\phi(g-f)$ is a square in $\mathcal{M}(Y,\mathcal{B})$. Hence $\phi(f)\leq \phi(g)$.
     Therefore any isomorphism between $\mathcal{M}(X,\mathcal{A})$ and $\mathcal{M}(Y,\mathcal{B})$ gives rise to an isomorphism of the positive cones $\mathcal{M}^+(X,\mathcal{A})$ and $\mathcal{M}^+(Y,\mathcal{B})$. 

    Conversely, if there is a semiring isomorphism $\phi\colon \mathcal{M}^+(X,\mathcal{A})\rightarrow \mathcal{M}^+(Y,\mathcal{B})$, then there is a homeomorphism $\psi\colon (\mathcal{RM}Cong(\mathcal{M}^+(X,\mathcal{A})),\{\mathfrak{m}^R(f,g)\colon f,g\in \mathcal{M}^+(X,\mathcal{A})\}) \rightarrow (\mathcal{RM}Cong(\mathcal{M}^+(Y,\mathcal{B})),\{\mathfrak{m}^R(h,k)\colon h,k\in \mathcal{M}^+(Y,\mathcal{B})\})$ defined as $\psi(\mathfrak{m}^R(f,g))=\mathfrak{m}^R(\phi(f),\phi(g))$. For realcompact measurable spaces $(X,\mathcal{A})$ and $(Y,\mathcal{B})$ we have, $\mathcal{RM}Cong(\mathcal{M}^+(X,\mathcal{A})\cong X\cong RMax(X,\mathcal{A})$ and $\mathcal{RM}Cong(\mathcal{M}^+(Y,\mathcal{B}))\cong Y\cong RMax(Y,\mathcal{B})$ (cf. Theorem \ref{P3, Th6.18}). Therefore $X\cong Y$ as measurable spaces.
 \end{proof}

 By \textbf{Meas}, we denote the category of measurable spaces and measurable functions. In particular, we are interested in the subcategory \textbf{RCTMeas}, consisting of realcompact $T$-measurable spaces. Theorem \ref{P3, Th6.19} solves the isomorphism problem for the semirings of the form $\mathcal{M}^+(X,\mathcal{A})$. We arrive at the following easy proposition. 

 \begin{proposition}
  
         The category \textbf{RCTMeas} is dual to the full subcategory of \textbf{CRig} (the category of commutative semirings and semiring homomorphisms) consisting of the semirings of the form $\mathcal{M}^+(X,\mathcal{A})$. 
         
 \end{proposition}
	
	\section*{Acknowledgement}
	The first author is grateful to the University Grants Commission (India) for providing a Junior Research Fellowship (ID: 211610013222/ Joint CSIR-UGC NET JUNE 2021).

	\bibliographystyle{plain}

\end{document}